\documentclass[12pt]{article}

\input xy
\xyoption{all}

\newtheorem{thm}{Theorem}[section]
\newtheorem{prop}[thm]{Proposition}
\newtheorem{lemma}[thm]{Lemma}
\newtheorem{cor}[thm]{Corollary}

\newcounter{ex}[section]

\newcommand{\Spec}{\mbox{$\mbox{Spec}$} }

\newcommand{\p}{\mbox{$p$}\ }

\newenvironment{eq}{\addtocounter{thm}{1}\begin{equation} }{\end{equation}}

\input xy
\xyoption{all}

\begin{document}
\title{Equivariant Euler characteristics and sheaf resolvents}

\author{Ph. Cassou-Nogu\`es, M.J. Taylor  }

\maketitle

\noindent{\bf Abstract:}
 For certain tame abelian covers of arithmetic surfaces $X/Y$ we obtain striking formulas, involving a quadratic form derived from intersection numbers, for the equivariant Euler characteristics of both the canonical sheaf $\omega_{X/Y}$ and also its square root $\omega_{X/Y}^{1/2}$. These formulas allow us us to carry out explicit calculations; in particular, we are able to exhibit examples where these two Euler characteristics and that of the the structure sheaf of $X$ are all different and non-trivial. Our results are obtained by using resolvent techniques together with the local Riemann-Roch approach developed in [CPT].
 
\medskip

\noindent {\bf \Large Introduction}

\bigskip
Let $N/K$ be a finite Galois extension of number fields with  Galois group $G$.  A number of interesting arithmetic modules
may be associated to such an extension: the ring of algebraic integers $\mathcal{O}_{N}$ of $N$; the codifferent $\mathcal{D}_{N/K}^{-1}$
of $N/K$; and, when the ramification subgroups of $N/K$ are of odd order, the square root of the codifferent
 $\mathcal{D}_{N/K}^{- {1}/{2}}$. Their structure as Galois modules has been studied  extensively.
 When  $N/K$ is at most tamely ramified, they are all three locally free ${\bf Z}[G]$-modules. It was proved in [T1] and [T2] that, for any tame Galois extension $N/K$, the classes of $\mathcal{O}_{N}$ and $\mathcal{D}_{N/K}^{-1}$ in  $ \mbox{Cl}({\bf Z}[G])$, the locally free class group of ${\bf Z}[G]$-modules, are  equal; that is to say, $\mathcal{O}_{N}$ is a selfdual ${\bf Z}[G]$-module.
 If in addition we suppose that  $G$ is of odd order, then in fact  $\mathcal{O}_{N}$,  $\mathcal{D}_{N/K}^{-1}$ and
 $\mathcal{D}_{N/K}^{-1/2}$ are all three free over $\mathbf{Z}[G]$.
   This result for the ring of integers is a consequence of the Fr\"ohlich conjecture, which was proved in [T3]. The result for the square root of the codifferent
  was obtained  in [ET]. In such a situation where one of the modules  $\mathcal{O}_{N}$,  $\mathcal{D}_{N/K}^{-1}$ and   $\mathcal{D}_{N/K}^{-1/2}$ are ${\bf Z}[G]$-free, we shall say that the module in question has a normal integral basis (abbreviated
  NIB).

  Over the past ten, or so, years a number of articles  have been  devoted to the study of the analogues of such Galois module problems in higher dimensions. The study of such questions for arithmetic surfaces is the central topic of this paper.

      Throughout this article $G$ denotes a finite group of $\it{odd}$ order . We consider a  $G$-cover $\pi: X\rightarrow Y$ of schemes which are projective and flat
  over $Spec(\mathbf{Z})$. For a  $G$-equivariant locally free sheaf $\mathcal{F}$ on $X$, the coherent cohomology
  groups $H^{i}(X, \mathcal{F})$ are finitely generated ${\bf Z}[G]$-modules which are of considerable interest. Their study leads us to consider the hypercohomology complex ${\bf R\Gamma }(X, \mathcal{F})$ of ${\bf Z}[G]$-modules. When
  the action of $G$ is tame, i.e. when for any point $x$ of $X$ the order $e_{x}$ of the inertia group $I_{x}$ of $x$
   is coprime to the residue field characteristic of $x$, this  complex  is perfect; that is to say it is isomorphic, in the derived category of complexes of ${\bf Z}[G]$-modules, to a bounded
   complex of projective ${\bf Z}[G]$-modules. The fact that this complex is perfect, under such hypotheses, was first shown in [CE] (see also [CEPT3]
   and  [P], Section 2). This then leads  to the following definition:
   \medskip

    \noindent {\it The sheaf $\mathcal{F}$ is said to have a normal integral basis when the complex
   ${\bf R\Gamma }(X, \mathcal{F})$,  in the derived category of complexes of ${\bf Z}[G]$-modules,  is isomorphic to a bounded
   complex of free
   ${\bf Z}[G]$-modules  }.
   \medskip

  \noindent Following Chinburg in [C1], we can associate to $\mathcal{F}$ an equivariant Euler characteristic $\chi^{P}( \mathcal{F})$
  in $\mbox{Cl}({\bf Z}[G])$, which measures the obstruction to the existence of a NIB for $\mathcal{F}$. More precisely, the sheaf $\mathcal{F}$ has a NIB
  if and only if $\chi^{P}(\mathcal{F})$ is trivial. We observe that in the situation described above where we take  $X= \mbox{Spec}(\mathcal{O}_{N})
  $,  it follows from the above mentioned theorems that the coherent sheaves associated to
   $\mathcal{O}_{N}$,  $\mathcal{D}_{N/K}^{-1}$ and   $\mathcal{D}_{N/K}^{-1/2}$ all have a NIB,  since their Euler
   characteristics are precisely the classes of these modules in  $\mbox {Cl}({\bf Z}[G])$.

   We  now introduce  the analogous sheaves of modules
   in the  geometric setting. We start by recalling that the different divisor of $X/Y$ is the divisor on $X$ given by
    $$D_{X/Y}=\sum_{x}(e_{x}-1)x$$
   with $x$ running over the set of codimension $1$ points of $X$. We shall be particularly  interested in the following $G$-equivariant sheaves: the structural sheaf $\mathcal{O}_{X}$; the canonical sheaf
   $\omega _{X/Y}=\mathcal{O}_{X}(D_{X/Y})$ of the cover $X\rightarrow Y$; and its square root
   $\omega _{X/Y}^{1/2}=\mathcal{O}_{X}(\frac {1}{2}D_{X/Y})$.

   The main object of our study is the existence, or otherwise, of a  NIB for these sheaves
   by  computing and  comparing  their  equivariant Euler
   characteristics.  In the case when $G$ is abelian and
    $X$ is a $G$-torsor  over $Y$,  the morphism   $\pi: X\rightarrow Y $
    is etale and hence the sheaves $\omega _{X/Y}$ and $\omega _{X/Y}^{1/2}$ both coincide with the structural sheaf $\mathcal{O}_{X}$, and furthermore they have a NIB by a theorem of Pappas in [P]. We shall see presently that when the cover $\pi: X\rightarrow Y $ is not etale then our sheaves do not necessarily have NIB.
    A new approach to such questions has been recently developed in [CPT]. Our results depend strongly on the use of this  paper, and illustrate how their new techniques permit the efficient calculation of such Euler characteristics.

     We now describe the contents of the paper. In Section 1 we describe the situation that we wish to study; here we introduce our notation and we state our main results.
     In  Section 2 we present some general results on equivariant duality; the content of this section derives from some working notes of T. Chinburg and G. Pappas, and we are most grateful for their permission to use their work here.
    In Section 3 we introduce the notions of a sheaf resolvent  and a divisor resolvent; these two concepts play a central role throughout this paper.
    The main theorems are proved in Sections 4 and 5.  We conclude with the detailed study of some examples in Section 6; we are extremely grateful to Arnaud Jehanne for his help with the computations in 6.c.

    \section {Notation and main results}

  Let $G$ be a finite abelian group of  odd order $n$, and  let $R$ denote either a Dedekind domain or  a complete discrete valuation ring; in all cases we denote the field of fractions of $R$ by
   $K$. In the case where $R$ is a valuation ring, we shall assume that $R$ contains the $n$-th roots of unity and that its  residue field $k$ is
  perfect and of characteristic prime to $n$. We consider a regular flat projective scheme $Y\rightarrow S=\mbox {Spec}(R)$. The fibres are of constant dimension and we denote this fibral dimension by  $d$. Let $\pi: X\rightarrow Y$ be a $G$-cover which is generically
  a $G$-torsor on $Y$ with  $X$ is regular.  Note that it follows from the assumptions that $\pi: X\rightarrow Y$ is flat (see Remark 3.1.a in [CPT]). When $R={\bf Z}$ we suppose that the ramification locus of this cover is supported on a finite set of rational primes  $\Sigma$ which is disjoint with the set of prime divisors of the order of $G$.

    In order to state our results we now suppose that $R={\bf Z}$ and $d=1$. We consider a $G$-equivariant,  coherent and invertible sheaf $\mathcal{F}$ on $X$.
    For any $p\in \Sigma $
   we denote by ${\bf Z}_{p}'$ the subring of ${\bf Q}_{p}^{c}$ obtained by adjoining the $n$-th
    roots of unity to ${\bf Z}_{p}$. By forming the base changes of  $\pi: X\rightarrow Y$  by
     $\mbox{Spec}({\bf Z}_{p}) \rightarrow \mbox{Spec}({\bf Z})$ and by $\mbox{Spec}({\bf Z}_{p}') \rightarrow \mbox{Spec}({\bf Z})$,  we obtain
     $G$-covers $\pi_{p}: X_{p}\rightarrow Y_{p}$ and $\pi_{p}': X_{p}'\rightarrow Y_{p}'$. We denote by  $\mathcal{F}'_{p}$
      the $X_{p}'$- sheaf  obtained from $\mathcal{F}$ by pullback.

      For an equivariant sheaf $\mathcal{F}'_{p}$ as above, for any  ${\bf Q}_{p}^{c}$-character $\varphi $ of $G$ and for a codimension one point $y$ of $Y'_p$,  in Section 2 we will define a rational number  $v_{y}( \mathcal{F}'_{p,\varphi })$, which
   depends on the ramification of $y$  in the cover
   $X_{p}'\rightarrow Y_{p}'$. We shall then use these rational numbers to define the
   {\it  local resolvent divisor  of $\mathcal{F}$ at $p$} by setting
   \begin{eq}
   r_{p}(\mathcal{F}, \varphi )=\sum_{y}v_{y}( \mathcal{F}'_{p,\varphi })y ,
   \end{eq}
   where $y$ runs over the set of codimension one points of $Y_{p}'$ which are contained in the special fiber $Y_{p}'^{(s)} $ of
   $Y_{p}' \rightarrow \mbox{Spec}({\bf Z}_{p}')$.  We may consider $r_{p}(\mathcal{F}, \varphi )$ as a vector  with rational coordinates.
    Presently we will see that
   $n v_{y}( \mathcal{F}'_{p,\varphi })$ is an integer for any such $y$, so that $n.r_{p}(\mathcal{F}, \varphi )$ is a divisor of
   $Y_{p}'$. These resolvent divisors play a similar role to that of Lagrange resolvents in the algebraic number field setting.
   For codimension one points $y$ and $z$  of $Y_{p}'^{(s)}$,  we denote their intersection number by  $y\cdot{z}$.  Recall that this integer is   the degree of the line
   bundle $\mathcal{O}_{Y'_{p}}(y)$ restricted to $z$. Let  $\omega _{Y_{p}'}$ be the canonical sheaf of
   $Y_{p}'\rightarrow {\bf Z}_{p}'$ and denote its first Chern class by $c_{1}(\omega _{Y_{p}'})$ . Then we define the integer $c_{1}(\omega _{Y_{p}'})\cdot{y}$ as the
   degree of the $0$-cycle $c_{1}(\omega _{Y_{p}'})\cap y$ of $Y_{p}'$ (see  Chapter 2 in [Fu]). Finally we set:
   \begin{eq}
   r_{p}(\mathcal{F}, \varphi )^{2}=\sum_{y, z}v_{y}( \mathcal{F}'_{p,\varphi })v_{z}( \mathcal{F}'_{p,\varphi })y\cdot{z} \ ,
   \end{eq}
   \begin{eq}
   c_{1}(\omega _{Y_{p}'})\cdot r_{p}(\mathcal{F}, \varphi )=\sum_{y}v_{y}( \mathcal{F}'_{p,\varphi })c_{1}(\omega _{Y_{p}'})\cdot{y},
   \end{eq}
   and
   \begin{eq}
   T_{p}(\mathcal{F}, \varphi)=r_{p}(\mathcal{F}, \varphi )^{2}+c_{1}(\omega _{Y_{p}'})r_{p}(\mathcal{F}, \varphi ).
   \end{eq}
   We observe that $r_{p}(\mathcal{F}, \varphi )^{2}$ may be thought of as the quadratic form,  defined  by  the intersection matrix,  evaluated on a
   local resolvent divisor; while $c_{1}(\omega _{Y_{p}'})\cdot r_{p}(\mathcal{F}, \varphi )$ may be thought of as a linear form,  evaluated on the same local resolvent divisor.

   Suppose now that there exists a locally free $\mathcal{O}_Y$-sheaf  $\omega _{Y/S}^{ {1}/{2}}$  with the property that
   $\omega _{Y/S}^{{1}/{2}}\otimes \omega _{Y/S}^{{1}/{2}}=\omega _{Y/S};$ we shall refer to $\omega _{Y/S}^{{1}/{2}}$ as a square root of $\omega _{Y/S}^{{1}/{2}}$.
    We then define a twist $\tilde \mathcal{F}$ of  $\mathcal{F}$  by setting
    \begin{eq}
     \tilde \mathcal{F}=\mathcal{F}\otimes \pi ^{*}(\omega _{Y/S}^{ {1}/{2}})  \ .
     \end{eq}
     It follows from the adjunction formula that $\tilde  \omega _{X/Y}^{ {1}/{2}}$ has the property that
      $\tilde\omega _{X/Y}^{{1}/{2}}\otimes \tilde\omega _{X/Y}^{{1}/{2}}=\omega _{X/S}$
    and we  therefore  denote  this sheaf by  $\omega _{X/S}^{ {1}/{2}}$.
    \medskip

     Once and for all we fix a sufficiently large finite Galois extension $E$ of ${\bf Q}$. We denote  the group
  of finite ideles of $E$ by $J_{f}(E)$; we let $R_{G}$ denote the additive group of the virtual $E$-characters of $G$ and we let $t$  denote the
  group homomorphism
  \begin{eq}
   t: \mbox{Hom}_{G_{\bf Q}}(R_{G}, J_{f}(E)) \rightarrow \mbox{Cl}({\bf Z}[G])
   \end{eq}
  of Fr\"ohlich's so-called $``$ Hom-description $"$,
  (see  Chapter 2 in [F1]).
    For each finite prime $l$ of ${\bf Z}$ we fix an embedding $j_{l}: {\bf Q}^{c} \rightarrow {\bf Q}_{l}^{c}$  and we denote  the closure of $j_{l}(E)$ in ${\bf Q}_{l}^{c}$ by
   $E_{l}$. The group $J_{l}(E)=(E\otimes {\bf Q}_{l})^{\times}$ may be identified with the Galois submodule of $J_{f}(E)$
   consisting of the finite ideles  of $E$ which are equal to $1$ outside $l$; we shall therefore view
    $\mbox{Hom}_{G_{\bf Q}}(R_{G}, J_{l}(E))$ as a subgroup of $\mbox{Hom}_{G_{\bf Q}}(R_{G}, J_{f}(E))$.
 We therefore obtain, via $j_{l}$,
   a homomorphism $j_{l}: E\otimes {\bf Q}_{l} \rightarrow  E_{l}$; we also obtain an isomorphism $\varphi \mapsto
   (\varphi )^{j_{l}}$ between the ${\bf Q}^{c}$-characters of $G$
    and the ${\bf Q}_{l}^{c}$-characters of $G$ and therefore an isomorphism between $R_{G}$ and $R_{G,l}$,  the additive groups of their virtual characters.  It is well known,
    (see for instance  Lemma II.2.1 of [F1]),  that the homomorphism

   $$j_{l}^{*}: \mbox{Hom}_{G_{\bf Q}}(R_{G}, J_{l}(E))\rightarrow \mbox{Hom}_{G_{{\bf Q}_{l}}}(R_{G,l}, E_{l}^{*})$$
   defined by $j_{l}^{*}(f)(\varphi )=f(\varphi ^{j_{l}^{-1}})^{j_{l}}$ is an isomorphism. We shall often abbreviate the notation $j_{l}^{*}(f)$
   to $f^{*}$.
    \noindent Our main results are the following:

    \begin{thm}\label{general} With the above notation and under the above assumptions on $X, Y$,  $G$,  and for both
    $\mathcal{F}=\omega _{X/Y}$ or $\omega _{X/Y}^{1/2}$ we have the following equality in the class group $\mbox{Cl}({\bf Z}[G])$:

    $$2\chi^{P} ( \mathcal{O}_{X})-2\chi^{P} (\mathcal{F})
    =\prod\limits_{\p\in \Sigma } t( g_p)$$
     where $ g_{p}^*$ is the  element of $\mbox{Hom}_{G_{{\bf Q}_{p}}}(R_{G,p}, E_{p}^{*})$ defined on irreducible
     ${\bf Q}_{p}^{c}$-characters $\varphi $ of $G$ by
     $$ g_{p}^{*}(\varphi )=p^{T_{p}(\mathcal{F}, \varphi)-T_{p}(\mathcal{O}_{X}, \varphi)} \ . $$
    \end{thm}

   \begin{thm} With the above notation and under the above assumptions on $X, Y$,  $G$, suppose there is an $\mathcal{O}_Y$-sheaf \  $\omega _{Y/S}^{1/2}$ as above; then for both
    $\mathcal{F}=\omega _{X/Y}$ or $\omega _{X/Y}^{1/2}$ we have the following equality in the class group $\mbox{Cl}({\bf Z}[G])$
    $$2\chi^{P} ( \tilde \mathcal{O}_{X})-2\chi^{P} (\tilde  \mathcal{F})
    =\prod\limits_{\p\in \Sigma} t(\tilde  g_p)$$
    where $ \tilde g_{p}^*$ is the  element of $\mbox{Hom}_{G_{{\bf Q}_{p}}}(R_{G,p}, E_{p}^{*})$ defined on irreducible ${\bf Q}_{p}^{c}$-characters $\varphi $ of $G$ by
     $$\tilde  g_{p}^{*}(\varphi )=p^{r_{p}(  \mathcal{F},  \varphi )^{2}
     -r_{p}( \mathcal{O}_{X}, \varphi )^{2}}\ \ . $$

   \end{thm}

\noindent {\bf Remarks.}

$\mathbf{1}.$ The factor $2$ in our formulas derives from the need to consider the determinant
of cohomology of bundles of even rank in [CPT]; note that this factor  can be removed in certain situations. (See Section 3 in [CPT] for further details.)

$\mathbf{2}.$ Since the sheaves $\mathcal{F}$, considered in the above theorems, are isomorphic to the structural sheaf
 $\mathcal{O}_{X}$  on the
  general fiber,  it should be possible to apply the techniques  developed in Section 8 of [CEPT1] to compute the difference
  between their Euler characteristics.  This would lead, under certain assumptions, to a generalisation of our theorems to  non-abelian $G$. The advantage of our method is that it yields an attractive formula, involving quadratic elements
  associated to resolvent divisors. One can hope that the two methods will provide quite different formulas, which it will be interesting to compare. This is something that we plan to do in the near future.

\medskip

  We denote by $\mathcal{M}$ the maximal order of ${\bf Q}[G]$ and we let $D({\bf Z}[G])$ be the kernel of
 the group homomorphism $\mbox{Cl}({\bf Z}[G])\rightarrow \mbox{Cl}(\mathcal{M})$ induced by extension of scalars.
 Using the representative for the class $2\chi(\mathcal{O}_{X})-2\chi(\omega _{X/Y})$ in $\mbox{Hom}_{G_{\bf Q}}(R_{G}, J_{f}(E))$
 given in Theorem 1.7 above, in Section 5 we will show:
  \begin {thm} Let $G$ be a finite abelian $l$-group.

\noindent i)   The order of the class $2\chi^{P}(\mathcal{O}_{X})-2\chi^{P}(\omega _{X/Y})$ is a power of $l$.

\noindent ii)  If $l$ is a regular prime number, then $2\chi^{P}(\mathcal{O}_{X})-2\chi^{P}(\omega _{X/Y})\in{D({\bf Z}[G]})$.

\noindent iii) If $G$ is of order $l$, then  $2\chi^{P}(\mathcal{O}_{X})=2\chi^{P}(\omega _{X/Y})$.
\end {thm}

\noindent We observe that this theorem provides a generalisation  for arithmetic surfaces (up to a factor of $2$) of Theorem 3.1 of [T{2}].
\medskip

\noindent {\bf Remarks. }

 \noindent {\bf 1.} It follows from Theorem 4.3.2  in [CPET1] that the results of this last  theorem remain true
if we replace $\chi^{P}(\omega _{X/Y})$ by $\chi^{P}(\Omega1 _{X/{\bf Z}})$.

\noindent {\bf 2.} The class group $\mbox{Cl}({\bf Z}[G])$ carries a natural duality involution $x\mapsto x^{D}$ (see Section 2).
The anti-selfduality of $\chi^{P}(\omega _{X/{\bf Z}}^{1/2})$ in  dimension $2$ will be shown in Section 2 to follow from
 a  general result proved for  arbitrary dimension (see Corollary 2.6).

\bigskip

 Let  $p\equiv  1$ mod  $ 24$ be a prime number and let $X_{1}(p)$ be the  model over $\mbox {Spec}({\bf Z})$ of the modular curve
       associated to the congruence subgroup $\Gamma _{1}(p)$ as described in  Section 4 in [CPT]. The group
        $\Gamma  =({\bf Z}/p{\bf Z})^{*}/\{\pm 1\}$ acts faithfully on $X_{1}(p)$. For a prime divisor $l$ of $p-1$, with $l>3$,
        we let  $H$  be the subgroup  of $\Gamma $ of index $l$ and   $G=\Gamma /H$ .
        We consider the projective schemes   $X=X_{1}(p)/H$ and $Y=X_{1}(p)/\Gamma $. Then
        $\pi: X\rightarrow Y$ is a $G$-cover to which we can  apply our general results.
        Let $\mathcal{P}$ be the prime ideal of $ {\bf Q}(\zeta _{l})$ defined by the chosen
       embedding $j_{p}:{\bf Q}^{c}\rightarrow{{\bf Q}_p^c}$ and we let $\bar \mathcal{P}$ denote its image under complex conjugation. For $1\leq a<l$ we denote by
       $\sigma_{a}$ the ${\bf Q}$-automorphism of ${\bf Q}(\zeta _{l})$ induced by $\zeta_{l}\mapsto \zeta_{l}^{a}$.
       Since $G$ is of prime order, by a theorem of Rim, any non-trivial abelian character  of $G$ induces an isomorphism between $\mbox{Cl}({\bf Z}[G])$ and
       $\mbox{Cl}({\bf Z}[\zeta _{l}])$. For $i\in \{1, 2\}$ we  define the elements $s_{i}$ in the group ring ${\bf Z}[G] $ by
       $$s_{i}=\sum_{1\leq a<l/2}a^{i}\sigma_{a}^{-1}.$$

\begin{thm}
        For a suitable choice of non-trivial abelian character $\varphi$ of $G$, we have equalities in  $\mbox{Cl}({\bf Z}[\zeta _{l}])$:

         \noindent i) $ \varphi (2\chi^{P}(\omega _{X/Y}^{1/2})-2\chi^{P}(\mathcal{O}_{X}))=
         [\mathcal{P}\overline \mathcal{P}]^{\frac{p-1}{12l}s_{1}} $.

         \noindent ii) $  \varphi (2l\chi^{P}(\omega _{X/Y}^{1/2}))
        =[\mathcal{P}\overline \mathcal{P}]^{\frac{p-1}{12l}s_{2}}$ and $
        \varphi (2l\chi^{P}(\mathcal{O}_{X}))
        =[\mathcal{P}\overline \mathcal{P}]^{\frac{p-1}{12l}(s_{2}-ls_{1})} $
       .
        \end{thm}
        Using this we immediately deduce:
        \begin {cor}

         Let $h_{l}^{+}$ be the class number of the maximal real subfield of ${\bf Q}(\zeta _{l})$ and
        assume that $h_{l}^{+}=1$; then we have:
        $$2\chi^{P} (\mathcal{O}_{X})=2\chi^{P} (\omega _{X/Y}^{1/2}) ,\ \ 2l\chi^{P} (\mathcal{O}_{X})=2l\chi^{P} (\omega _{X/Y}^{1/2})=0\ .
        $$
        If in addition  $l$ is a regular prime number, then
        $$2\chi (\mathcal{O}_{X})=2\chi (\omega _{X/Y})=2\chi(\omega _{X/Y}^{1/2})=0.$$
        \end{cor}

        \noindent {\bf Remark. }
        1. The Lefschetz-Riemann-Roch theorems of [CEPT2] and the techniques developed in [P] provide efficient tools for determining the prime-to-$l$ part of the Euler characteristic of the sheaves that we have considered here.

        \smallskip

        This corollary can be used to provide  families of covers $X\rightarrow Y$ for which our sheaves have Euler characteristic of order two . In  Section 6 we will see that
        Theorem 1.10 can  also be used  to construct
          families of examples where such Euler characteristics have order greater than 2. For instance for $p=182857$ and $l=401$ the cover $X\rightarrow Y$ provides us with a tame cover
          such that: $\chi^{P} (\mathcal{O}_{X}), \ \chi^{P}(\omega _{X/Y})$ and $\chi^{P} (\omega _{X/Y}^{1/2})$  all have order greater than 2; where
          $2\chi^{P} (\mathcal{O}_{X})=2\chi^{P}(\omega _{X/Y})$; but where $\chi^{P} (\mathcal{O}_{X})\not\neq \chi^{P}(\omega _{X/Y}^{1/2})$.

\section {An equivariant Duality Theorem.   }

      \noindent In this section we do not impose any restriction of the dimension of $Y$ and we no longer suppose $G$ to be abelian. However, we note that if $G$ is abelian, then there is a simpler proof of the results of this section by following [RD] and working over $\mbox{Spec}(\mathbf{Z}[G])$ as the base.

      The action of complex conjugation on the characters of $G$ induces an involutary
      automorphism on $\mbox{Cl}({\bf Z}[G])$.  Thus  if $[M]$ is an element of $\mbox{Cl}({\bf Z}[G])$ represented by
      $f\in \mbox{Hom}_{G_{\bf Q}}(R_{G}, J(E))$,  we define $\overline {[M]}$ as the element of  $\mbox{Cl}({\bf Z}[G])$ represented by $
      \overline {f}$ where for $\chi \in R_{G}$
      $$\overline f(\chi)=f(\overline \chi) \ $$
      and one checks that this automorphism maps $ker(t)$ in 1.6 into itself. We denote by $\mbox{Cl}({\bf Z}[G])^{+}$ (resp. $\mbox{Cl}({\bf Z}[G])^{-}$) the subgroup of elements of $x \in
      \mbox{Cl}({\bf Z}[G])$ with the property that $\bar x=x$ (resp. $\bar x=-x$).

      Clearly the group automorphism $f\mapsto f^{-1}$ on $ \mbox{Hom}_{G_{\bf Q}}(R_{G}, J(E))$ induces an involution   on
      $\mbox{Cl}({\bf Z}[G])$. By composing these two involutions we obtain a further involution which we denote by $[M]\mapsto [M]^{D}$.
      In Proposition 3 of Appendix A, IX in [F1],  Fr\"ohlich gives an interpretation of this latter involution by proving  that
      for any locally free ${\bf Z}[G]$-module $M$ one has the equality:
      \begin{eq}
      [M]^{D}=[M^{D}]
      \end{eq}
      where $M^{D}=\mbox{Hom}_{{\bf Z}}(M, {\bf Z})$ is the $\bf Z$-linear dual of $M$.
      Such involutions were used to provide some of the first restrictions on the Galois module structure of algebraic rings of integers and their ideals. In the previous section we associated  to  certain $G$-equivariant sheaves $\mathcal{F} $ on $X$ an equivariant Euler
      characteristic $\chi^{P}( \mathcal{F})$ in $\mbox{Cl}({\bf Z}[G])$;  it is then natural to consider the image of this class
      under the above involution. The initial  aim of this section is to give an interpretation of the class
      $\chi^{P}(\mathcal{F})^{D}$ in terms of a duality functor for $Y$-sheaves, which extends Fr\"ohlich's above result. We provide just such a description in Corollary 2.5, and we note that the relative
      dimension $d$ of $X$ over ${\mbox{Spec}(\bf Z) }$ appears explicitly in this result. The result that we give is derived from  an equivariant generalization   of the duality theorem for projective morphisms
      given in Theorem 11.1 in  III of [RD]. This  $``$equivariant$"$ version  was first given in some unpublished notes of Ted Chinburg and George Pappas.

        We start by
       introducing  a small amount of further notation. Our main references are Section 2 of [C2]  and Section 2 of [P]. As previously,  $G$ is a finite
      group and $ Y$ is a projective flat scheme over ${\bf Z}$.  Let $K(Y, G)$ (resp. $K({\bf Z}, G)$)
      be the homotopy category of $\mathcal{O}_{Y}[G]$-modules (resp. ${\bf Z}[G]$-modules)  and let
      $K^{+ }(Y, G)$ (resp. $K^{+}({\bf Z}, G)$ be the subcategory of complexes in $K(Y, G)$ (resp. $K({\bf Z}, G)$) which are bounded
      below and which have coherent (resp. finitely generated) cohomology. We let $D(Y, G),\ D^{+}(Y, G),\ D({\bf Z}, G)$ and
      $D^{+}({\bf Z}, G)$ be their respective  derived categories . We define $K(Y),\ K^{+}(Y),\ D(Y)$ and $D^{+}(Y)$
      (resp. $K({\bf Z}),\ K^{+}({\bf Z}),\ D({\bf Z})$
      and $D^{+}({\bf Z})$) similarly by considering complexes of $\mathcal{O}_{Y}$-modules (resp. ${\bf Z}$-modules ).
      The global section functor $\Gamma $ has a right derived functor
      $$R\Gamma ^{+}: D^{+}(Y, G)\rightarrow D^{+}({\bf Z}, G)\ $$
      (see Section 2 in II of [RD]).

      In Theorem 1.1 of [C1], Chinburg proved:
      \begin{thm}
      If $F^{\bullet}$ is a bounded complex in $K^{+}(Y, G)$ with the property that each
      stalk of each term of $F^{\bullet}$ is a cohomologically trivial $G$-module,  then $R\Gamma ^{+}(F^{\bullet})$ is isomorphic in
      $D^{+}({\bf Z}, G)$ to a finite complex  of  finitely generated projective ${\bf Z}[G]$-modules.
      \end{thm}

      This fact led him to associate to any such complex an  Euler characteristic $\chi^P(R\Gamma ^{+}(F^{\bullet}))$ in the
      class group $\mbox{Cl}({\bf Z}[G])$.
      When $\pi: X\rightarrow Y$ is a tame $G$-cover, then, under the hypotheses of Section 1,  one can prove  that  if
      $\mathcal{F}$ is any $G$-equivariant coherent sheaf on $X$, then  the complex of $D^{+}(Y, G)$, which is $\pi_{*}(\mathcal{F})$
      in degree $0$ and $0$ elsewhere, satisfies the conditions required by Chinburg's theorem. In this case the class
      $\chi^{P}( \mathcal{F})$, referred to in the Introduction,  coincides with the class $\chi^{P}(R\Gamma ^{+}(\pi_{*}(\mathcal{F})))$
       defined above.
       \medskip

      \noindent For any $F^{\bullet}$ in $K(Y, G)$ and any $T^{\bullet}$ in $K(Y)$ we have the complex
      ${\bf Hom}^{\bullet}_{\mathcal{O}_{Y}}(F^{\bullet}, T^{\bullet})$ in $K(Y, G)$   (see Section 3 in II of[RD]).
      We thereby  obtain a bifunctor:
      $${\bf Hom}^{\bullet}_{\mathcal{O}_{Y}}: K(Y, G)^{0}\times K(Y)\rightarrow K(Y, G)\  $$
      and a derived bifunctor
      $$R{\bf Hom}^{\bullet}_{\mathcal{O}_{Y}}: D(Y, G)^{0}\times D^{+}(Y)\rightarrow D^{+}(Y, G) \  $$
      where, as usual, the superscript $0$ denotes the opposite category. The bifunctors ${\bf Hom}^{\bullet}_{{\bf Z}}$ and
       $R{\bf Hom}^{\bullet}_{{\bf Z}}$ are defined in the same way. For any $F^{\bullet}$ in $D(Y, G)$ and $T^{\bullet }$
      in $D^{+}(Y)$ we may then consider the Euler characteristic
      $\chi^P(R\Gamma ^{+}(R{\bf Hom}^{\bullet}_{\mathcal{O}_{Y}})(F^{\bullet}, T^{\bullet}))$
      in $\mbox {Cl}({\bf Z}[G])$. The following theorem is a special case of the above mentioned equivariant  duality theorem:

      \begin {thm} (Chinburg, Pappas.) Let $F^{\bullet}$ be a bounded complex in $K^{+}(Y, G)$ as in Theorem 2.2. Let
      $h: Y\rightarrow S= Spec ({\bf Z})$ be the structural morphism. Then one has the equality of Euler
      characteristics in $Cl({\bf Z}[G])$

      $$\chi^P(R\Gamma ^{+}(R{\bf Hom}^{\bullet}_{\mathcal{O}_{Y}}(F^{\bullet}, h^{!}(\mathcal{O}_{S}))))=
      \chi^P(R{\bf Hom}^{\bullet}_{{\bf Z}}(R\Gamma ^{+}(F^{\bullet}), {\bf Z}))\ .$$
      \end{thm}

      \noindent {\bf Remark} The construction of $h^{!}$ is given in III of [RD] .
      \medskip

    To give a very brief idea of the proof of the theorem we remark that the equality of these two Euler characteristics is a consequence
    of the existence of an isomorphism $\Gamma (\theta_{h})$ in $D^{+}({\bf Z})$:
    $$\Gamma (\theta_{h}): R\Gamma ^{+}(R{\bf Hom}^{\bullet}_{\mathcal{O}_{Y}}(F^{\bullet}, h^{!}(\mathcal{O}_{S})))
    \rightarrow  R{\bf Hom}^{\bullet}_{{\bf Z}}(R\Gamma ^{+}(F^{\bullet}), {\bf Z})\ $$
    and this latter isomorphism is  obtained by applying the global section functor to the duality isomorphism $\theta_{h}$ of Section 11 in III of [RD].

    \begin {cor} Since $Y$ is flat over $S$, the fibres all have constant dimension which we denote by $d$. We then have the equality:

     $$(-1)^{d}\chi^P(R\Gamma ^{+}(R{\bf Hom}^{\bullet}_{\mathcal{O}_{Y}}(F^{\bullet}, \omega _{Y/S})))=
      \chi^P(R{\bf Hom}^{\bullet}_{{\bf Z}}(R\Gamma ^{+}(F^{\bullet}), {\bf Z}))\ .$$
      \end {cor}

      \begin{Proof}  It follows from [RD] that $h^{!}(\mathcal{O}_{S})$ is the complex $\omega _{Y/S}[d]$ whose only non-zero
      term is $\omega _{Y/S}$ in degree $-d$. Moreover we have the equality:
      $$
      \chi^P(R\Gamma ^{+}(R{\bf Hom}^{\bullet}_{\mathcal{O}_{Y}}(F^{\bullet}, \omega _{Y/S}[d])))=
      (-1)^{d} \chi^P(R\Gamma ^{+}(R{\bf Hom}^{\bullet}_{\mathcal{O}_{Y}}(F^{\bullet}, \omega _{Y/S})))\ . $$
      The corollary therefore follows immediately from the theorem.
      \end{Proof}

      \begin {cor} With the hypotheses and notation of the previous corollary, we have the following
      equality in $\mbox{Cl}({\bf Z}[G])$:

      $$\chi^{P}(\mathcal{F})^{D}=
      (-1)^{d}\chi^{P}(R\Gamma ^{+}(R{\bf Hom}^{\bullet}_{\mathcal{O}_{Y}}(\pi_{*}(\mathcal{F}), \omega _{Y/S})))\ . $$
      \end {cor}

      \begin{Proof} We consider the complex $\pi_{*}(\mathcal{F})$ whose only non-zero term is $\pi_{*}(\mathcal{F})$ in degree $0$.
      Since the action $(X, G)$ is tame we can apply Theorem 2.2 to deduce that the complex $R\Gamma ^{+}(\pi_{*}(\mathcal{F}))$ is perfect.
      Let $M^{\bullet}$ be a bounded complex of finitely generated projective ${\bf Z}[G]$-modules which is isomorphic to
      $R\Gamma ^{+}(\pi_{*}(\mathcal{F}))$ in $D^{+}({\bf Z}, G)$.
      Then $R{\bf Hom}^{\bullet}_{{\bf Z}}(M^{\bullet}, {\bf Z})={\bf Hom}(M^{\bullet}, {\bf Z})$ is
      the class in $D^{+}({\bf Z}, G)$ of the complex whose $(-j)^{th }$ term is $\mbox{Hom}_{{\bf Z}}(M^{j}, {\bf Z})=(M^{j})^{D}$.
      Therefore it follows from Fr\"ohlich's duality result that we have the equality in $\mbox{Cl}({\bf Z}[G])$:
      $$ \chi^{P}(R{\bf Hom}^{\bullet}_{{\bf Z}}(R\Gamma ^{+}(\pi_{*}(\mathcal{F})), {\bf Z}))=
      \chi^{P}( \mathcal{F})^{D}\ . $$
      The result now follows from the previous corollary.
      \end{Proof}
      \medskip

      \noindent We now wish use the above result to derive some properties for the Euler characteristics of the $G$-equivariant
      sheaves on $X$ that we have considered.

      \begin {cor} Suppose that  the hypotheses of the previous corollaries are satisfied. Then
      \item i) $$\chi^{P}(\mathcal{O}_{X})^{D}=(-1)^{d}\chi^{P}(\omega _{X/S}).$$

      \item ii) Moreover,  if  the ramification indices of $X\rightarrow Y$ are odd and if there exists an $\mathcal O_Y$-line bundle $\omega _{Y/S}^{1/2}$, with the property that $\omega _{Y/S}^{1/2}\otimes\omega _{Y/S}^{1/2}=\omega _{Y/S}$, then
      $$\chi^{P}( \omega _{X/S}^{1/2})^{D}=(-1)^{d}\chi^{P}( \omega _{X/S}^{1/2})\ . $$
      In particular, when $d=1$ the class  $\chi^{P}( \mathcal{O}_{X})-\chi^{P}( \omega _{X/S})$ belongs to $\mbox{Cl}({\bf Z}[G])^{-}$ and the class $\chi^{P}(\omega _{X/S}^{1/2})$ belongs
      to $\mbox{Cl}({\bf Z}[G])^{+}$.

       \end {cor}

       \begin{Proof} From the duality formula of Corollary 2.5 we see that, in order to prove this corollary, we need to evaluate
       $R{\bf Hom}^{\bullet}_{\mathcal{O}_{Y}}(\pi_{*}(\mathcal{F}), \omega _{Y/S})$ for
       $\mathcal{F}=\mathcal{O}_{X}$ and
       $\mathcal{F}=\omega _{X/S}^{1/2}$. Since the action  $(X, G)$ is tame,   in both cases $\pi_{*}(\mathcal{F})$ is a locally free
       $\mathcal{O}_{Y}[G]$-module.
       Therefore $R{\bf Hom}^{\bullet}_{\mathcal{O}_{Y}}(\pi_{*}(\mathcal{F}), \omega _{Y/S})=
       {\bf Hom}_{\mathcal{O}_{Y}}(\pi_{*}(\mathcal{F}), \omega _{Y/S})$, and so now
       we consider ${\bf Hom}_{\mathcal{O}_{Y}}(\pi_{*}(\mathcal{F}), \omega _{Y/S})$. First we observe that
       $${\bf Hom}_{\mathcal{O}_{Y}}(\pi_{*}(\mathcal{F}), \omega _{Y/S})\simeq
      {\bf Hom}_{\mathcal{O}_{Y}}(\pi_{*}(\mathcal{F}), \mathcal{O}_{Y})\otimes_{\mathcal{O}_{Y}}\omega _{Y/S}. $$
      From Propositions 4.25. and 4.32 in VI of [L]  we deduce  that when $\mathcal{F}=\mathcal{O}_{X}$ we have
      $${\bf Hom}_{\mathcal{O}_{Y}}(\pi_{*}(\mathcal{O}_{X}), \mathcal{O}_{Y})\simeq \pi_{*}(\omega _{X/Y})\ . $$
      Therefore we have isomorphisms of $\mathcal O_{Y}[G]$-modules
      $${\bf Hom}_{\mathcal{O}_{Y}}(\pi_{*}(\mathcal{O}_X), \omega _{Y/S})\simeq
      \pi_{*}(\omega _{X/Y})\otimes _{\mathcal{O}_{Y}}\omega _{Y/S}\simeq \pi_{*}(\omega_{X/S} ) $$
      with the latter isomorphism following from the adjunction and the projection formulas. This then shows $i)$.

      Now let
      $\mathcal{F}=\omega _{X/S}^{1/2}$. From the projection formula we deduce  that
      $\pi_{*}(\omega _{X/S}^{1/2})=\pi_{*}(\omega _{X/Y}^{1/2})\otimes_{\mathcal{O}_{Y}}\omega _{Y/S}^{1/2}$. Therefore
      $${\bf Hom}_{\mathcal{O}_{Y}}(\pi_{*}(\omega _{X/S}^{1/2}), \mathcal{O}_{Y})\simeq
     {\bf Hom}_{\mathcal{O}_{Y}}(\pi_{*}(\omega _{X/Y}^{1/2}), \mathcal{O}_{Y})
     \otimes _{\mathcal{O}_{Y}}{\bf Hom}_{\mathcal{O}_{Y}}(\omega _{Y/S}^{1/2}, \mathcal{O}_{Y}) \ . $$
      By use of the trace pairing we know that $\pi_{*}(\omega _{X/Y}^{1/2})$ is a self-dual $\mathcal O_{Y}[G]$-module; also,  since $\omega_{Y/S} ^{1/2}$ is invertible, we know that
     ${\bf Hom}_{\mathcal{O}_{Y}}(\omega _{Y/S}^{1/2}, \mathcal{O}_{Y})=\omega _{Y/S}^{-1/2}$.
     In conclusion we deduce from the above that
     $${\bf Hom}_{\mathcal{O}_{Y}}(\pi_{*}(\omega _{X/S}^{1/2}), \omega _{Y/S})
     \simeq \pi_{*}(\omega _{X/Y}^{1/2})\otimes _{\mathcal{O}_{Y}}\omega _{Y/S}^{-1/2}
     \otimes _{\mathcal{O}_{Y}}\omega _{Y/S}\simeq \pi_{*}(\omega _{X/S}^{1/2}) \ .$$
     This then completes the proof of the corollary.
     \end{Proof}

\section{Sheaf and divisor resolvents}

 In this section we assume  that $G$ is a finite group of exponent $n$ and that   $R$  is  a complete discrete valuation ring whose residue class field $k$ is of characteristic $p$, which is coprime  to $n$, and that $R$ contains the $n$th roots of unity. We assume that
   $Y\rightarrow S=\mbox{Spec}(R)$ is   of absolute dimension $d+1$. We  associate to  $G$ the  constant group scheme over $S$, which we again denote by $G_S$; we let
$G^{D}_{S}=\mbox{Spec}(R[G])$ be the Cartier dual of
 $G_S$, and $G^{D}_{Y}$ denotes the fiber product $G^{D}_{S}\times _{S}Y$.
   For any  scheme $Z$, we denote by  $\mathcal{K}_{Z}$ the sheaf of stalks of meromorphic functions on $Z$. This is
  the sheaf associated to the presheaf defined on affine open subschemes $\mathcal{U}$ of $Z$ by  $
 Frac (\mathcal{O}_{Z}(\mathcal{U}))$, where we let  $Frac(A)$ denote the ring of fractions of the ring $A$.

  Let $\mathcal{F}$ be a coherent invertible subsheaf  of $\mathcal{K}_{X}$ whose support is contained in the branch locus of the cover $\pi:X\rightarrow{Y}$. We assume that $D=\sum _{x} d_{x} x$ is a divisor of $X$ with $gD=D$ for any $g\in{G}$ so that $d_{x}$ depends only on the $G$-orbit of $x$ and furthermore
   $d_{x}=0$ whenever the inertia subgroup $I_{x}$ of $x$ is trivial. For future reference we note that by duality the canonical sheaf $\omega_{X/Y}$ has this property; moreover, if the ramification of $X/Y$ is odd, then by ? we know that  $\omega_{X/Y}^{1/2}$ also has this property.

   We consider the sheaf $\pi_{*}(\mathcal{F})$:
this is a $G$-equivariant, locally free $\mathcal{O}_{Y}$-module which may be viewed as an invertible
$G^{D}_{Y}$-module (see  for instance 5.17 in II of [H]).
Any character  $\varphi  : G \rightarrow  \Gamma (S, \mathcal{O}_{S}^{*})$
 of $G$  provides us with an
$S$-point of $G^{D}_{S}$ and hence,   by base change,  with   a $Y$-point  $\varphi: Y \rightarrow G^{D}_{Y}$
of $G^{D}_{Y}$.

\vskip 0.1 truecm

\noindent {\bf Definition 3.1. } We define the  resolvent sheaf associated to $\mathcal{F}$ and $\varphi$ to be the invertible $Y$-sheaf:
 $$\mathcal{F}_{\varphi }=\varphi ^{*}((\pi_{*})(\mathcal{F})) \ .$$

 \noindent $\mathcal{F}_{\varphi }$  is easily seen to be the subsheaf of $\pi _{*}(\mathcal{F})$ consisting
of those local sections on which $G$ acts by the character $\varphi $.

$n$-fold multiplication induces a homomorphism of $Y$-sheaves
$$\mu:
\pi_{*}(\mathcal{K}_{X})^{\otimes {n}} \rightarrow \mathcal{K}_{X}$$
which, by   restriction,   identifies
 $\mathcal{F}_{\varphi }^{\otimes n}$ with an invertible subsheaf of $\mathcal{K}_{Y}$ that may be defined by a divisor.
 The principal aim of  this section is to give an explicit expression for this divisor.

  Let $x$ be a codimension one point of $X$, let $I_{x}$  be  the inertia group  of $x$ and
  $e_{x}$ be the order of $I_{x}$. The action of $I_{x}$ on the cotangent space at $x$ defines a faithful character $\psi  _{x}$
 with values in $k^{*}$. Since $n$ is coprime to the characteristic of $k$ and since $R$ is a complete discrete valuation ring, we may view $\psi_x$ as taking values in $R^{*}$.
 Hence for any character $\varphi $ of $G$ there exists
a unique integer $n(\varphi , x), \ 0\leq n(\varphi , x)<e_{x}$, such that the  restriction to $I_{x}$ of
 $\varphi $  is equal to $\psi _{x}^{n(\varphi , x)}$. We set $\pi(x)=y$. Since the points of $x$ above $y$ are all conjugate
 under the action of $G$, both the group $I_{x}$ and the integer $n(\varphi , x)$ depend only on $y$. Moreover,
 as  noted previously,   since the divisor $D=\sum_{x}d_{x}x$ is $G$-invariant,  the integer
$d_{x}$ also depends only  on   $y$.  We shall therefore denote these objects  $I_{y}, n(\varphi ,y)$ and  $d_{y}$.
For any rational number  $a $ we denote its integral part by $[a]$  and its  fractional part by $\{a\}=a-[a]$.

\noindent We set
\begin{eq}
f_{\mathcal{F}}(\varphi , y)=\frac{d_{y}}{e_{y}}-\{ \frac {n(\varphi ,y)+d_{y}}{e_{y}}\} \ .
\end{eq}
\begin{prop} For any abelian character $\varphi $  of $G$ the map $\mu $ identifies  $\mathcal{F}_{\varphi }^{\otimes n}$ with
$\mathcal{O}_{Y}(F_{\mathcal{F}}(\varphi ))$ with
$$F_{\mathcal{F}}(\varphi )=\sum _{y} n f_{\mathcal{F}}(\varphi , y)y$$
where $y$ runs over the set of codimension one points of $Y$ which are contained in the special fiber of $Y \rightarrow S$.

\end{prop}

\begin{Proof} We  may  assume, for the purposes of the proof,  that the decomposition subgroup of $x$ is equal to the inertia subgroup.
We consider  a codimension one point $y$ on the special fiber of $Y$ and we fix a point $x$ on $X$
  such that $\pi(x)=y$. For the sake of simplicity we will  write $e$ for $e_{y}$, $d$ for
$d_{y}$, $I$ for $I_y$ etc.. Since $R$ is a complete discrete valuation ring with residue characteristic $p$ coprime to $n$ and which contains the $n$th roots of unity, we are in a tame Kummer situation and so we can choose a uniformising parameter   $\omega _{x}$   of $\mathcal{O}_{X, x}$
with the property that $\omega _{y}=\omega _{x}^{e}$ is a uniformiser of $\mathcal{O}_{Y, y}$ . From the very definition of
$\mathcal{F}$ we know that
$\mathcal{F}_{x}=\omega _{x}^{-d}\mathcal{O}_{X,x}$. We define $q$ and $r$  by the equality
$-d=qe+r$ subject to the restriction that  $0\leq r<e $. Then   $\omega _{x}^{-d}=\omega_{y}^{q}\omega_{x}^{r}$. Again, since we are in a tame Kummer situation,
 we know that  $\alpha =\frac {1}{e}(1+\omega _{x}+...\omega _{x}^{e-1})$ is
a free basis of $\mathcal{O}_{X, x}$ as an $\mathcal{O}_{Y, y}[I]$-module. It then follows easily
that $\omega_{y}^{q}\omega_{x} ^{r}\alpha $ is
a free basis of $\mathcal{F}_{x}$ over $\mathcal{O}_{Y, y}[I]$.
  This implies that  the stalk of $\mathcal{F}_{\varphi }^{n}$ at $y$ is given  by:
$$ \omega _{y}^{-n.f_{\mathcal{F}}(\varphi )}
\mathcal{O}_{Y,y}= \mathcal{F}_{\varphi, y }^{n}= (\omega_{y}^{q}\omega_{x} ^{r}\alpha \mid \psi _{x}^{n(\varphi )})^{n}\mathcal{O}_{Y,y} $$
where for $a\in \mathcal{F}_{x}$ and  a character $\theta $ of $I$,  $(a\mid \theta)$  denotes the Lagrange resolvent
$$(a\mid \theta)=\sum_{g\in I}a^{g}\theta(g^{-1}).$$
 By the standard properties of Lagrange resolvents and the definition of the cotangent character $\psi_x $ we have the following
 equalities:
$$(\omega_{y}^{q}\omega_{x} ^{r}\alpha \mid \psi _{x}^{n(\varphi )})=
\omega _{y}^{q}(\omega _{x}^{r}\alpha\mid \psi _{x}^{n(\varphi )})=
\omega_{y}^{q}\omega _{x}^{r}(\alpha \mid \psi _{x}^{n(\varphi )-r}) \ .$$
 By an easy
computation we obtain that
$$(\alpha \mid \psi_{x} ^{n(\varphi )-r})=\sum_{g\in{I_{y}}}\alpha^{g^{-1}}\psi_{x}^{n(\varphi)-r}(g)=$$
$$\frac {1}{e}\sum_{g\in I_{y}}\sum _{0\leq k< e}\omega_{x} ^{k}\psi_{x} ^{-k-r+n(\varphi )}(g)=
 \frac {1}{e}\sum_{0\leq k< e}\omega_{x} ^{k}\sum _{g\in I_{y}}\psi_{x} ^{-k-r+n(\varphi )}(g) \ .$$
 We conclude that $(\alpha \mid \psi_{x} ^{n(\varphi )-r})=\omega_{x} ^{n(\varphi )-r}$ if $r\leq n(\varphi )$
  and  $e+n(\varphi )-r$  otherwise. Piecing this together we obtain that
  $-n.f_{\mathcal{F}}(\varphi )= \frac {n}{e}(n(\varphi ) +eq)$ if
   $r\leq n(\varphi )$ and $\frac {n}{e}(n(\varphi ) +e(q+1))$ otherwise. The proposition now follows at once from  this  equality.
   \end{Proof}

 \noindent As we have seen in the proof of the above proposition, the rational numbers  $f_{\mathcal{F}}(\varphi ,y) $ are defined by the equalities
 \begin{eq}
  f_{\mathcal{F}}(\varphi ,y)=-\frac {1}{n}v_{y}(\mathcal{F}_{ \varphi, y }^{ n}) \ .
  \end{eq}
  From now on we will denote  this number by  $-v_{y}(\mathcal{F}_{\varphi })$.
  \vskip 0.1 truecm
  \noindent {\bf Definition 3.2.}  The resolvent divisor of $\mathcal{F}$ at $\varphi $ is defined to be
  \begin{eq}
   r(\mathcal{F}, \varphi )=\sum_{y}v_{y}( \mathcal{F}_{\varphi,y })y\
   \end{eq}
   where $y$ runs over the set of codimension one points of $Y$ which are contained in the special fiber of $Y \rightarrow S$.

   \noindent We note that $n.r(\mathcal{F}, \varphi )$ is a divisor on $Y$ with the property that
 $$\mathcal{F}_{\varphi }^{\otimes n}
   =\mathcal{O}_{Y}(-n.r(\mathcal{F},  \varphi))  \ .$$
    We will sometimes abuse notation and write   $\mathcal{F}_{\varphi }=\mathcal{O}_{Y}(-r(\mathcal{F},  \varphi ))$.
    \medskip

    \noindent {\bf Examples. } The following three examples are central to our study.
    \smallskip

    \noindent {\bf 1.}   If $\mathcal{F}$ is the structural sheaf of $X$, then $D=0$ and  from the above proposition it follows that
    \begin{eq}
  v_{y}(\mathcal{O}_{X, \varphi })=\frac{n(\varphi , y)}{e_{y}}\ \  \mbox{and}\  \ r(\mathcal{O}_{X},  \varphi )=\sum_{y}\frac{n(\varphi , y)}{e_{y}}y\ .
  \end{eq}
  Note that this is Lemma 3.5 of [CPT].
\smallskip

\noindent {\bf 2.}  We let $\mathcal{S}(\varphi)$
    denote the set of codimension one points $y$  of $Y$ contained in the special fiber of $Y \rightarrow S$,
  with the property that $n(\varphi, y)> 0$. We set
  \begin{eq}
  f(\varphi)=\sum_{y\in \mathcal{S}(\varphi)}y .
  \end{eq}
  If we take  $\mathcal{F}=\omega_{X/Y}$, then $d_y=e_y-1$ and we get
  $$r(\omega _{X/Y}, \varphi)=r(\mathcal{O}_{X}, \varphi)-f(\varphi) $$
  by using Proposition 3.2 and the fact that
  $$\{ \frac{n(\varphi,y)+e_y-1}{e_y} \}=\frac{e_y-1}{e_y}+\frac{n(\varphi,y)}{e_y}-1.$$

  \smallskip

  \noindent {\bf 3.}  We let
    $\mathcal{S}'(\varphi)$ denote the set of codimension one points $y$  of $Y$ contained in the special fiber of $Y \rightarrow S$,
  with the property that that
   $n(\varphi, y)>e_{y}/2$.  We set
  \begin{eq}
    f'(\varphi)=\sum_{y\in \mathcal{S}'(\varphi)}y \ .
  \end{eq}
  If we take  $\mathcal{F}=\omega^{1/2}_{X/Y}$, then $d_y=(e_y-1)/2$ and we get
  $$r(\omega^{1/2} _{X/Y}, \varphi)=r(\mathcal{O}_{X}, \varphi)-f'(\varphi) $$
  by using Proposition 3.2 and the fact that
  $$\{ \frac{n(\varphi,y)+e_(y-1)/2}{e_y} \}=\frac{e_y-1}{2e_y}+\frac{n(\varphi,y)}{e_y} \quad \mbox{if} \quad n(\varphi,y)<e_y/2,$$
   $$\{ \frac{n(\varphi,y)+e_(y-1)/2}{e_y} \}=\frac{e_y-1}{2e_y}+\frac{n(\varphi,y)}{e_y}-1 \quad \mbox{if} \quad n(\varphi,y)>e_y/2.$$

  From Proposition 3.2 it follows that
  $$r(\omega _{X/Y}^{1/2}, \varphi)=r(\mathcal{O}_{X}, \varphi)-f'(\varphi).$$

   \begin{cor}  Let $\varphi $ be an abelian character of $G$ and let $\bar \varphi $ be its complex conjugate. Then:

 i)  $$r(\mathcal{O}_{X}, \varphi )+r(\mathcal{O}_{X},  \bar \varphi )=f(\varphi) $$

  ii) $$r(\omega _{X/Y}^{1/2},  \varphi )+r(\omega _{X/Y}^{1/2}, \bar \varphi )=0$$

  iii) $$r(\omega _{X/Y}^{{1}/{2}},  \varphi )=r(\mathcal{O}_{X},  \varphi^{2})-
 r(\mathcal{O}_{X}, \varphi) \ . $$
 \end{cor}
 \begin{Proof}
 Part $(i)$ follows from Example 1 above and part $(ii)$ comes from Example 3. To prove part $(iii)$ we note that if $n(\varphi,y)<e_y/2$ then $n(\varphi2,y)=2n(\varphi,y)$ and if $n(\varphi,y)>e_y/2$ then $n(\varphi2,y)=2n(\varphi,y)-e_y$, and so
 $$r(\mathcal{O}_{X},  \varphi^{2})-r(\mathcal{O}_{X},  \varphi)=(\sum_y\frac{2n(\varphi,y)}{e_y})-f'(\varphi)$$
 $$=2r(\mathcal{O}_{X},  \varphi)-f'(\varphi)=r(\mathcal{O}_{X},  \varphi)+r(\omega _{X/Y}^{{1}/{2}},  \varphi ).$$
 \end{Proof}
 \medskip

 \noindent {\bf Remark.} We denote the  sheaf $\mathcal{O}_{Y}(-f(\varphi))$ by $F(\varphi)$.
 From the above corollary we deduce the following equalities of invertible $\mathcal{O}_Y$-sheaves:
 $$\mathcal{O}_{X, \varphi}^{n}\otimes \mathcal{O}_{X, \bar \varphi}^{n}= F(\varphi)^{n},
 \ \ ({\omega _{X/Y, \varphi}^{1/2}})^{n}\otimes ({\omega _{X/Y, \bar \varphi}^{1/2}})^{n}=\mathcal{O}_{Y}\ $$
 and we observe that the first equality provides us with a geometric analogue of  Theorem 18 in [F{2}].

\section {Euler characteristics and representatives.   }

 \subsection {The Riemann-Roch Theorem} Let    $R$ be  a Dedekind domain. Henceforth we suppose that  $ Y\rightarrow S=\mbox{Spec}({R})$ is of absolute
 dimension $2$ and we consider a $G$-cover of $Y$ satisfying the assumptions of Section 1.

  To any  locally free coherent $G$-sheaf  $\mathcal{F}$ on $X$ we can associate the following two invariants: the first is the equivariant Euler characteristic
$\chi^{P} (\mathcal{F})\in\mbox {Cl}({R}[G])$, which was introduced in Section 1 (see [C1] for further details).  We now briefly recall the
construction of the second invariant that we wish to study. Recall that the sheaf $\pi_{*}(\mathcal{F})$ may be viewed as a locally free coherent $G_{Y}^{D}$-sheaf.
Let $\tilde h: G_{Y}^{D}\rightarrow G_{S}^{D}$ be the base change of $h:Y\rightarrow S$ by $G_{S}^{D}\rightarrow S$.
Then the total derived  image $R\tilde h(\pi_{*}(\mathcal{F}))$ in the derived category of complexes of $G_{S}^{D}$-sheaves is represented by a perfect complex.
Therefore,  following Knudsen and Mumford,  we can define the invertible ${R}[G]$-module  $\mbox{det}(R\tilde h(\pi_{*}(\mathcal{F}))$ and
consider its class  in $\mbox {Pic}(G_{S}^{D})=\mbox{Pic}({R}[G])$. Since $G$ is abelian  this latter group may  be
identified with the class group $\mbox{Cl}({R}[G])$ and we have the following equality in the class group $\mbox {Cl}({R}[G])$, (see ?2.c  in [P]  and also Section 3 of [CPT]):
$$\chi^{P}( \mathcal{F})=[\mbox{det}(R\tilde h(\pi_{*}(\mathcal{F}))]\ .$$
In fact  we shall be interested in twice this class and so, with the notation of [CPT], we set
$$\delta (\mathcal{F})=\mbox{det}(R\tilde h(\pi_{*}(\mathcal{F}))^{\otimes 2}.$$

\noindent For $i, 1\leq i\leq 2,$ we let $\mathcal{F}_{i}$ denote a coherent,  invertible,  $G$-equivariant sheaf on $X$,  defined by a
divisor $D_{i}$. We assume that $D_{1}\leq D_{2}$ and that both these divisors are supported on the branch locus of $\pi: X\rightarrow Y$.
Our aim is to describe
$\psi (\mathcal{F}_{1}, \mathcal{F}_{2})=
2(\chi^{P}( \mathcal{F}_{2})-\chi^{P}( \mathcal{F}_{1}))$ in $\mbox{Cl}(R[G])$.  It follows from the above that
$\psi (\mathcal{F}_{1}, \mathcal{F}_{2})=\delta (\mathcal{F}_{2})\otimes \delta (\mathcal{F}_{1})^{-1}$.
Since $D_{1}\leq D_{2}$ we have a natural injection of sheaves $f: \pi_{*}(\mathcal{F}_{1})\rightarrow \pi_{*}(\mathcal{F}_{2})$
and so we have the equality
$$\psi (\mathcal{F}_{1}, \mathcal{F}_{2})=[\delta (\pi_{*}(\mathcal{F}_{1})\rightarrow \pi_{*}(\mathcal{F}_{2}))]$$
where $\pi_{*}(\mathcal{F}_{1})\stackrel{f}{\rightarrow }\pi_{*}(\mathcal{F}_{2})$ denotes the perfect complex
of $G_{Y}^{D}$-sheaves with terms in positions $-1$ and $0$. We shall denote by $\delta (f)$ the locally free $R[G]$-module
defined by the ``$\delta $-image'' of this complex; that is to say, the determinant of the  push down, via $f$,  of this complex in $\mbox{Cl}(R[G])$. With this notation, the previous equality may be re-written as
$$\psi (\mathcal{F}_{1}, \mathcal{F}_{2})=[\delta(f)]\ .$$ Observe that $\delta(f)$ actually defines a submodule of $K[G]$. Since
$G$ is abelian, the determination of this module reduces to the computation of $\varphi (\delta(f))$ for each abelian character
$\varphi$ of $G$.  By pullback, any such character induces a morphism of sheaf resolvents
$$\varphi^{*}(f): \mathcal{F}_{1,\varphi}\rightarrow \mathcal{F}_{2, \varphi}\ $$ and, by the functorial
properties of the determinant of the cohomology,  it follows that we have $\varphi(\delta(f))=\delta (\varphi^{*}(f))$.
\medskip

\noindent The module $\delta (f)$ may of course be described by its local components, and so we now consider the local situation where $R$ is a complete discrete valuation ring.
 We introduce a small amount of further notation. Let $k$ denote the residue class field of $R$.  For any scheme  $f: T\rightarrow \mbox{Spec}(k)$  we let  $G_{0}(T)$ denote
  the Grothendieck group of coherent sheaves on $T$. We identify $G_{0}(\mbox{Spec}(k))$ with $\bf Z$. Assuming $f$ to be proper,  we
   define the Euler characteristic of a coherent sheaf $\mathcal{F}$ on $T$ by the equality:
  $$\mathcal{X}(T, \mathcal{F})=f_{*}([\mathcal{F}]) \  $$
  where $[\mathcal{F}]$ is the class of  $\mathcal{F}$ in  $G_{0}(T)$.

  \noindent We let $h_{s}: Y_{s}\rightarrow \mbox{Spec}(k)$ denote the special fiber of $Y$. The complex
$$c_{\varphi}: \mathcal{F}_{1,\varphi}\stackrel {\varphi^{*}(f)}\rightarrow \mathcal{F}_{2, \varphi}$$ is then a complex of locally free
$\mathcal{O}_Y$-modules which is exact off $Y_{s}$. Thus the cokernel of ${\varphi^{*}(f)}$,  which we denote by
$\mathcal{T}_{\varphi}(\mathcal{F}_{1}, \mathcal{F}_{2})$,  is a coherent sheaf of $\mathcal{O}_Y$-modules which is supported on $Y_{s}$.
It therefore defines a class in the Grothendieck group $G_{0}(Y_{s})$ of coherent $Y_{s}$-modules. Let $\lambda$ be a uniformizing parameter of
$R$. Then it is easily shown that
$$\varphi (\delta (f))= R\lambda  ^{-2\mathcal{X}(\mathcal{T}_{\varphi}(\mathcal{F}_{1}, \mathcal{F}_{2}))}\ $$
where we view the Euler characteristic $\mathcal{X}(\mathcal{T}_{\varphi}(\mathcal{F}_{1}, \mathcal{F}_{2}))\in \mathbf Z$ as above. The remainder
of this section will be devoted to the computation of the Euler characteristic
$\mathcal{X}(\mathcal{T}_{\varphi}(\mathcal{F}_{1}, \mathcal{F}_{2}))$ for  particular  choices of $\mathcal{F}_{1}$ and
$\mathcal{F}_{2}$.
\medskip

   For any integer $m\geq 0$ we let
  $A_{m}(T)$ denote the group of $m$-cycles on $T$ modulo rational equivalence and we set
  $$A_{*}(T)=\oplus_{0\leq m}A_{m}(T)\ .$$
  We write $A_{*}(T)_{\bf{Q}}=A_{*}(T)\otimes{\bf{Q}}$ and,  as usual, we identify $A_{0}(\mbox{Spec}(k))$ with $\bf Z$. The general Riemann-Roch theorem provides us with homomorphisms
  $\tau _{T}:G_{0}(T)\rightarrow A_{*}(T)_{\bf Q}$ which are covariant for proper morphisms. For any coherent sheaf
   $\mathcal{F}$ on $T$, the Riemann-Roch theorem states that
   $$\mathcal{X}(T, \mathcal{F})=f_{*}(\tau _{T}([\mathcal{F}]))\ .$$
   In this equality the right-hand side is the push forward of the zero component of $\tau _{T}([\mathcal{F}])$ to
   $A_{0}(Spec(k))_{{\bf Q}}$ identified with ${\bf Q}$.
   \smallskip

    For any invertible $\mathcal{O}_Y[G]$-sheaf $\mathcal{F}$ and any abelian ${\bf Q}_{p}^{c}$-character $\varphi$ of $G$, as previously,  we
   consider the  resolvent  divisor  $r(\mathcal{F}, \varphi)=\sum_{y}v_{y}(\mathcal{F}_{\varphi})y$,  (see  Definition 3.2 in Section 3), and
   we define $T(\mathcal{F}, \varphi )$ to be the integer
   $$T(\mathcal{F}, \varphi )=r(\mathcal{F}, \varphi)^{2}+ c_{1}(\omega _{Y/S})\cdot r(\mathcal{F}, \varphi)\ . $$
    Recall that for two codimension one points  $y$, $z$ on  the special fiber $Y_{s}$ of $Y$,  we denote their
    intersection number by  $y\cdot{z}$  and   $c_{1}(\omega _{Y/S})\cdot{y}$ is the
   degree of the $0$-cycle $c_{1}(\omega _{Y/S})\cap y$ on $Y$ (see Chapter 2 of[Fu] ). In the case
   where we have a ``square root'' $\omega _{Y/S}^{1/2}$ of $\omega _{Y/S}$, in 1.5 we have defined the twist $\tilde \mathcal{F}$ of $\mathcal{F}$ by the rule
   $$\tilde \mathcal{F}=\mathcal{F}\otimes \pi^{*}(\omega _{Y/S}^{1/2}).$$

  \begin{prop} Let $\mathcal{F}$ denote either $\omega _{X/Y}$ or $\omega _{X/Y}^{1/2}$. Then for any abelian
  ${\bf Q}_{p}^{c}$-character $\varphi $ of $G$  one has the following equalities:

 \noindent  i)
 $$2\chi(\mathcal{T}_{\varphi }(\mathcal{O}_{X}, \mathcal{F}))=
 T(\mathcal{F}, \varphi)-T(\mathcal{O}_{X}, \varphi)\ .$$

  \noindent ii) If there exists  a ``square root'' of $\omega _{Y/S}$, then
 $$ 2\chi(\mathcal{T}_{\varphi }( \tilde \mathcal{O}_{X}, \tilde \mathcal{F}))=r(\mathcal{F}, \varphi )^{2}
 -r(\mathcal{O}_{X },  \varphi )^{2}   \ .$$
   \end{prop}

 \noindent Proof.   Since the proofs of these equalities are similar,  we shall only  give the proof of  $ii)$ when $\mathcal{F}=\omega _{X/Y}$
 and will leave the proof of the remaining cases
 to the reader.
  Let $\varphi $ be an abelian character of $G$ and let $c_{\varphi}$ be the complex  of $Y$-line bundles
  $$\mathcal{O}_{X, \varphi }\otimes_{\mathcal{O}_{Y}}\omega _{Y/S}^{1/2}\stackrel {\varphi^{*}(\tilde f)}\rightarrow
  \omega _{X/Y, \varphi}\otimes_{\mathcal{O}_{Y}}\omega _{Y/S}^{1/2} $$ at degrees $-1$ and $0$.
  Since $c_{\varphi}$ is exact outside $Y_{s}$ and
  provides  a locally free resolution of
  $ \mathcal{T}_{\varphi }(\tilde \mathcal{O}_{X}, \tilde \mathcal{\omega }_{X/Y})$,
  (henceforth abbreviated to $\tilde \mathcal{T}_{\varphi}$),
   it follows from the Riemann-Roch Theorem (see Theorems 18.2 and 18.3 in [Fu])
  that

  $$\chi (\tilde \mathcal{T}_{\varphi} )=(h_{s})_{*}((ch_{Y_{s}}^{Y}(c_{\varphi})\cap Td(h))_{0})$$
   where $(h_{s})_{*}$ is the push forward
  $A_{0}(Y_{s})_{\bf Q}\rightarrow A_{0}(\mbox{Spec}(k))_{\bf Q}=\bf Q$, $ch_{Y_{s}}^{Y}(c_{\varphi})$ is the localized Chern character
  which lives in the bivariant group $A(Y_{s}\rightarrow Y)_{\bf Q}$, and $Td(h)$ is the Todd class of $h$,  which is an element of
  $A_{*}(Y)_{\bf Q}$ that we will describe presently.  We first observe   that the complex $c_{\varphi}$ may be written as
   $\omega _{Y/S}^{1/2}\otimes d_{\varphi}$,
  where $d_{\varphi}$ is the complex with terms in degrees $-1$ and $0$.
  $$\mathcal{O}_{X, \varphi }\stackrel {\varphi^{*}( f)}\rightarrow \omega _{X/Y, \varphi} \ .$$
   It therefore follows from  Proposition 18.1 (c) in [Fu] that
    $ch_{Y_{s}}^{Y}(c_{\varphi})=ch(\omega _{Y}^{1/2})\cap ch_{Y_{s}}^{Y}(d_{\varphi})$. We now consider the complex $d_{\varphi}^{\otimes n}$.
      For the sake of simplicity we  set
    $F(\varphi )=-n.r(\mathcal{O}_{X},  \varphi )$ and $R(\varphi )=n.f(\varphi)$. It follows from Proposition 3.2 that $d_{\varphi}^{\otimes n}$ is the complex
    $\mathcal{O}_{Y}(F(\varphi) ) \rightarrow \mathcal{O}_{Y}(F(\varphi )+R(\varphi ))$. We  observe
    that once again  $d_{\varphi}^{\otimes n}$ can be written as the tensor product $\mathcal{O}_{Y}(F(\varphi ))\otimes b_{\varphi}$ where $b_{\varphi}$ is the new
    complex $\mathcal{O}_{Y}\rightarrow  \mathcal{O}_{Y}(R(\varphi ))$.  Since $R(\varphi )$ is an effective divisor and since we
     assume that $X$ and $Y$ are relative curves, it follows  from Propositi
on 3.10 (b) and (e) in [CPT]  that
    $$ch_{Y_{s}}^{Y}(b_{\varphi})=  [R(\varphi )]+\frac {[R(\varphi )]^{2}}{2} $$ and therefore
    $$ch_{Y_{s}}^{Y}(d_{\varphi}^{\otimes n})=[R(\varphi )]+ c_{1}(\mathcal{O}_{Y}(F(\varphi )))\cap [R(\varphi )]
    + \frac {[R(\varphi )]^{2}}{2}.$$ Since for any integer $q$ we know that
     $ch_{Y_{s}}^{Y, q}(d_{\varphi}^{\otimes n})=n^{q}ch_{Y_{s}}^{Y, q}(d_{\varphi})$,
    it follows that
    $$ch_{Y_{s}}^{Y}(d_{\varphi})=\frac {[R(\varphi )]}{n}+ \frac {c_{1}(\mathcal{O}_{Y}(F(\varphi )))\cap [R(\varphi )]}
    {n^{2}}+ \frac {[R(\varphi )]^{2}}{2n^{2}} \ .$$
    Since $c_{\varphi}=\omega _{Y/S}^{1/2}\otimes d_{\varphi}$,  we finally obtain that
    $$ch_{Y_{s}}^{Y}(c_{\varphi})= \frac {[R(\varphi )]}{n}+ \frac {[R(\varphi )] \cap c_{1}(\omega _{Y/S}^{1/2})}{n}
    + \frac {c_{1}(\mathcal{O}_{Y}(F(\varphi )))\cap [R(\varphi )]}
    {n^{2}}+ \frac {[R(\varphi )]^{2}}{2n^{2}} \ .$$
    From the very definition of $Td(h)$ we deduce that

    $$Td_{1}(h)=-\frac {c_{1}(\omega _{Y/S})}{2}=-c_{1}(\omega _{Y/S}^{1/2})$$ (see (3.e) in [CPT]).
    Therefore, piecing  the above together, we obtain that
    $$2\chi (\tilde \mathcal{T}_{\varphi })=2(ch_{Y_{s}}^{Y}(c_{\varphi})\cap Td(h))_{0}=
    2\frac {c_{1}(\mathcal{O}_{Y}(F(\varphi )))\cap [R(\varphi )]}
    {n^{2}}+ \frac {[R(\varphi )]^{2}}{n^{2}} \ .$$
    Since we have the equalities
    $$\frac {[R(\varphi )]^{2}}{n^{2}}=f(\varphi)^{2} \ \ \mbox{and}\ \  \frac {c_{1}(\mathcal{O}_{Y}(F(\varphi )))\cap [R(\varphi )]}
    {n^{2}}=-r(\mathcal{O}_{X}, \varphi )\cdot f(\varphi)\ , $$ it suffices to use the equalities of
    Example 2 in Section 3 in order to deduce the required formula from above .

    \subsection {Proof of Theorems 1.7 and 1.8}
    Our hypotheses and our notations are  those of Theorems 1.7 and 1.8.  Since the proofs are similar we will give the proof of
   Theorem 1.7 and leave the proof of Theorem 1.8 to the reader. Our aim is to find  a representative for the class
    $  2(\chi^{P}(  \mathcal{F})-\chi^{P}(\mathcal{O}_{X}))$; that is to say an element $g$
     in the group $\mbox {Hom}_{G_{\bf Q}}(R_{G}, J_{f}(E))$ with the property that
     $$2(\chi^{P}(  \mathcal{F})-\chi^{P}(\mathcal{O}_{X}))=t(g)\ .$$
     For any element $p$   of $\Sigma $ we check that the map
     $$u_{p}: \varphi \mapsto p^{T_{p}(\mathcal{O}_{X}, \varphi)-T_{p}(\mathcal{F}, \varphi)}$$
    belongs to   $\mbox{Hom}_{G_{{\bf Q}_{p}}}(R_{G,p}, E_{p}^{*})$. Therefore we can define  $g_{p}$ to be  the  unique element of
     $ \mbox{Hom}_{G_{{\bf Q}}}(R_{G}, J_{p}(E))$ such that $g_{p}^{*}=u_{p}$. We set $g=\prod_{p\in \Sigma}g_{p}$.
    Our aim now is  to show that $g$ is a representative for our class.

    Since we have a group isomorphism
 $${\bf Q}_{p}[G]^{*}\simeq \mbox{Hom}_{G_{{\bf Q}_{p}}}(R_{G,p}, E_{p}^{*})$$
     given by  $x\mapsto (\varphi \mapsto \varphi (x))$, we may define
     $a_{p}$ to be the unique element of ${\bf Q}_{p}[G]^{*}$ corresponding to
     $u_{p}$ under the above isomorphism. By taking $a_{p}=1$ for $p\notin \Sigma $ and $a_{p}$ as above for $p\in \Sigma$, we obtain
     a finite idele $(a_{p})_{p}$ of ${\bf Q}[G]$ and thus an element in $\mbox{Pic}({\bf Z}[G])$ by considering the class of the fractional ideal
     $\cap_{p}{\bf Z}_{p}[G]a_{p}\cap {\bf Q}[G]$. This class corresponds to  $t(g)$ under the identification between $\mbox{Cl}({\bf Z}[G])$
     and $\mbox{Pic}({\bf Z}[G])$. Under this identification, as seen previously,  $2(\chi^{P}(  \mathcal{F})-\chi^{P}(\mathcal{O}_{X}))$
     identifies with $\psi (\mathcal{O}_{X}, \mathcal{F})$. Therefore, in order to prove the theorem, it suffices to prove that
     for each finite place $p$ we have the equality $\psi (\mathcal{O}_{X}, \mathcal{F}){\bf Z}_{p}={\bf Z}_{p}[G]a_{p}$. When
     $p\notin \Sigma$ this follows from the definition of $\psi (\mathcal{O}_{X}, \mathcal{F})$. For $p\in \Sigma$,
     since $p$ is coprime to $n$, it suffices to show that,  for any abelian ${\bf Q}_{p}^{c}$-character $\varphi$ of $G$,  we have
     $\varphi (\psi (\mathcal{O}_{X}, \mathcal{F}){\bf Z}'_{p})=\varphi(a_{p}){\bf Z}'_{p}$, where,  as previously, ${\bf Z}'_{p}$ is obtained
     by adjoining the $n$-th roots of unity to ${\bf Z}_{p}$. Since the functor $\delta $ commutes with base change, this last equality
     is equivalent  to the equality:
     $$\psi (\mathcal{O}_{X'_{p}}, \mathcal{F}'_{p})=
     p^{{T_{p}(\mathcal{O}_{X}, \varphi)-T_{p}(\mathcal{F}, \varphi)}}{\bf Z}'_{p}
     =p^{T(\mathcal{O}_{X'_{p}}, \varphi)-T(\mathcal{F}'_{p}, \varphi)}{\bf Z}'_{p}\ .
     $$
     This now  follows from Proposition 4.1.(i) because $p$ is a uniformizing parameter of ${\bf Z}'_{p}$.

      \section {The class $\chi(\mathcal{O}_{X})-\chi(\omega _{X/Y})$  }

        The aim of this section is to prove Theorem 1.9.   We now assume that
       $G$ is an abelian group of order $l^{N}$ where $l$ is a prime number, $l\geq 3$ and   $l\notin \Sigma$.  For any $p\in \Sigma$ we let ${\bf Z}'_{p}$ be
       the discrete valuation ring obtained by adjoining to ${\bf Z}_{p}$ a primitive $l^{N}$-th root of unity  and
       we let $\mathbf{Q} _p(\zeta_{l^N})$ denote its  field of fractions. Let
               $\pi_{p}: X'_{p}=X\otimes_{\bf Z}{\bf Z}'_{p}\rightarrow Y'_{p}=Y\otimes_{\bf Z}{\bf Z}'_{p}$
       denote the $G$-cover obtained from $\pi : X\rightarrow Y$
       by base change. There is then an action of  $V_{p}=\mbox{Gal}(K_{p}/{\bf Q}_{p})$ on ${\bf Z}'_{p}$,  and hence on $X'_{p}$ and $Y'_{p}$,
       and therefore on the codimension one points of $X'_{p}$ and $Y'_{p}$. Since the actions of $G$ and $V_{p}$ on $X'_{p}$ commute,
       it follows that for any $\omega \in V_{p}$ and any codimension one  point $x$ of
    $X$,   one has $I_{x^{\omega}}=I_{x}$ and $\psi _{x^{\omega}}=\psi _{x}^{\omega}$ (recall that  $\psi _{x}$ is the
    character of $I_{x}$ defined by the action of $I_{x}$ on the cotangent space in $x$). This implies that for any
    abelian character $\varphi: G\rightarrow {\bf Z}^{'\times}_{p}$,  $\omega \in V_{p}$
     we will have  $n(\varphi^{\omega}, x^{\omega})=n(\varphi, x)$. As previously, when $x\mapsto y$ we shall write $I_{y}$ and
    $n(\varphi, y)$ for $I_{x}$ and $n(\varphi, x)$.

     Let  $E$ be the  number field obtained by adjoining
 to $\bf Q$ the $l^{N}$-th roots of unity.
In Theorem 1.7 we proved  that
$$2(\chi(\mathcal{O}_{X})-\chi (\omega _{X/Y}))=\prod_{p\in \Sigma}t(g_{p}) $$
where $g_{p}$ is the element of $\mbox{Hom}_{G_{\bf Q}}(R_G, J_{p} (E))$ defined by
$$g_{p}^{*}(\varphi)=p^{T_{p}(\omega _{X/Y},\varphi)-T_{p}(\mathcal{O}_{X}, \varphi)}.\ $$

Let  $p$ be an element of $ \Sigma $ and let  $S(\varphi)$ denote the set of codimension one points $y$ of $Y'_{p}$ contained in the special fiber of $Y'_{p}\rightarrow Spec({\bf Z}'_{p})$ such that $n(\varphi, y)> 0$.
We know from Example 2 of Section 3  that, for any abelian character $\varphi$ of $G$
$$r_{p}(\omega _{X/Y}, \varphi)=r_{p}(\mathcal{O}_{X}, \varphi)
-f(\varphi)\ \  . $$
 Moreover, in 3.e. of [CPT] we find
 that $$c_{1}(\omega _{Y _{p}})\cdot{y}=-y^{2}-2\chi(y, \mathcal{O}_{y}).$$ Using 1.4 and Example 2 in Section3, we can use the above  to deduce that
$$T_{p}(\omega _{X/Y},\varphi)-T_{p}(\mathcal{O}_{X}, \varphi)=
f(\varphi)^{2}+\sum_{y\in S(\varphi)}(y^{2}+2\chi(y, \mathcal{O}_{y}))-2f(\varphi)\cdot r_{p}(\mathcal{O}_{X}, \varphi)$$
where,  as before, $y\cdot{z}$ (resp. $y^{2}$) denotes the intersection  (resp. self-intersection) number of $y$ and $z$ (resp. $y$).  We set:
$$ a(\varphi)=f(\varphi)^{2}+\sum_{y\in S(\varphi)}(y^{2}+
2\chi(y, \mathcal{O}_{y}))\ .$$
Since $a(\varphi)-2f(\varphi)r(\mathcal{O}_{X}, \varphi)$ is an Euler characteristic, and since by definition $a(\varphi)$ is an integer,
it follows that $2f(\varphi)\cdot r_{p}(\mathcal{O}_{X}, \varphi)$ must also be an integer. We observe that  for any
element $\omega \in \mbox{Gal} ( {\bf Q}_{p}^{c}/{\bf Q}_{p})$ there exists $r_{\omega}$, coprime to $l$, such that $\varphi^{\omega}=\varphi^{r_{\omega}}$.
We therefore deduce that $S(\varphi)=S(\varphi^{\omega})$ and hence $a(\varphi)=a(\varphi^{\omega})$. This implies that
$g_{p}$ may be written in $\mbox{Hom}_{G_{\bf Q }}(R_{G}, J_{p}(E))$ as a product $u_{p}h_{p}$ where for any abelian character $\varphi$
we define $u_{p}^{\ast}(\varphi)=p^{\ a(\varphi)}$ and $h_{p}^{\ast}(\varphi)=p^{-2f(\varphi)\cdot r(\mathcal{O}_{X_{p}}, \varphi)}$.

The theorem will be a consequence of the following  two  lemmas.
\begin {lemma} $t(u_{p})$ belongs to $D({\bf Z}[G])$.
\end {lemma}

\begin{Proof} The key observation is that for any abelian character $\varphi$ of $G$ and any $\omega \in G_{\bf Q}$ from the above we know that
$$a(\varphi_{p})=a((\varphi^{\omega})_{p})\ .$$
It therefore follows that the character map  $v_p$, whose value on an abelian character $\varphi$ is given by $v_{p} (\varphi) = p^{\ a(\varphi_{p})}$, belongs to $\mbox{Hom}_{G_{\bf Q}}(R_{G}, \bf Q^{\times})$.
Let $w_{p}$ be the element of $\mbox{Hom}_{G_{\bf Q}}(R_{G}, J(E))$ defined on abelian characters of $G$ by the rule:
 $$  w_{p}(\varphi)_{v} =  \left\{ \begin{array}{ccc}
1  \  \mbox{if }\ v\  \mbox {divides }\  p\\

p^{-a(\varphi_{p})}\ \mbox{if } v \ \mbox{if not }
\end{array} \right. $$
and so by definition  $w_{p} \in  \mbox{Hom}_{G_{\bf Q}}(R_{G}, U(E))$. Since
$(v_{p}w_{p})^{\ast}=u_{p}^\ast$ we deduce that $v_{p}w_{p}=u_{p}$. Therefore
$$t(u_{p})=t(v_{p})t(w_{p})=t(w_{p})\ .$$
Since $t(w_{p})\in D({\bf Z}[G])$ the lemma is proved.
\end{Proof}
\medskip

\noindent If $F$ is a number field and $x$ is a finite idele of $F$, then the content of $x$ is a fractional ideal of $F$, which we denote by $c(x)$.

\begin {lemma}

i)   There exists an integer $m \geq 0$ such that $l^{m}t(h_{p})$ belongs to $D({\bf Z}[G])$.

ii) If $G$ is of order $l$, then $t(h_{p})$ belongs to $D({\bf Z}[G])$.
\end{lemma}

\begin{Proof} In order to prove $ i) $, it suffices to show  that for any non-trivial abelian character $\varphi$ of $G$ there exists an
integer $s\geq 0$ such that $l^{s}c(h_{p}(\varphi))$ is a principal ideal.
  For such a character  $\varphi $, let $F$ be the field ${\bf Q}(\varphi)$ and let $\mathcal{P}$ be the prime ideal of $F$ defined
  by the restriction of $j_{p}$
to $F$. It
follows from the definition that $v_{\mathcal{P}}(c(h_{p}(\varphi)))=-2f(\varphi_{p})r_{p}(\mathcal{O}_{X}, \varphi_{p})$.
The Galois group $H_{p}$ of the extension ${\bf Q}_{p}(\varphi_{p})/{\bf Q}_{p}$ can be identified via  $j_{p}$ with the decomposition group of $\mathcal{P}$.
For any integer $a$,  coprime to $l$,   we denote by $\sigma_{a}$ the automorphism of $E$ defined by the property that  on a primitive $l^{N }$-th root of
unity $\zeta $ we have
$\sigma_a(\zeta)=\zeta^{a}$.
We then let  $\{ a_{i}, 1\leq i\leq q\}$ be a  set
of integers  such that the $\{\sigma_{a_{i}}\}$ are a set of representatives of $\mbox {Gal}(F/{\bf Q})/H_{p}$.
 We have the equalities
$$v_{\mathcal{P}^{\sigma_{a}^{-1}}}(c(h_{p}(\varphi)))=v_{\mathcal{P}}((c(h_{p}(\varphi))^{\sigma_{a}})
=v_{\mathcal{P}}(c(h_{p}(\varphi^{a}))=-2f(\varphi_{p}^{a})\cdot r_{p}(\mathcal{O}_{X}, \varphi_{p}^{a})\ . $$
Since it is clear  that  $f(\varphi_{p}^{a})=f(\varphi_{p})$, we deduce  from Example 1 in Section 3 that
$$f(\varphi_{p}^{a})\cdot{r_{p}(\mathcal{O}_{X}, \varphi_{p}^{a})}=
\sum _{y\in S(\varphi_{p})}(f(\varphi_{p})\cdot y){\frac{n(\varphi_{p}^{a}, y)}{e_y}}\ .$$
Therefore we obtain
$$c(h_{p}(\varphi))=\mathcal{P}^{\alpha(\varphi)}\ \mbox{with}\
 \alpha(\varphi)=-2\sum_{y\in S(\varphi_{p})}(f(\varphi_{p})\cdot y)\sum_{1\leq i\leq q}\frac{n(\varphi_{p}^{a_{i}}, y)}{e_{y}}
\sigma_{a_{i}}^{-1}\ .$$
The group $V_{p}=Gal(\mathbf{Q} _p(\zeta_{l^N})/{\bf{Q}}_{p})$ acts on $S(\varphi_{p})$.   We let $\{y_{j}, 1\leq j\leq r\}$ be a set of orbit representatives of
$S(\varphi_{p})$ under the action of $V_{p}$. For any $j, \ 1\leq j\leq r$, let $V_{p,j}$ denote the isotropy group of the codimension one  point $y_{j}$. Since
$V_{p,j}$ is a subgroup of ${\bf Q}_{p}(\zeta_{l^{N}})/{\bf Q}_{p}(\psi _{y_{j}})$, we note that it must be an $l$-group.
 For any $u\in V_{p}$ and any $y\in S(\varphi_{p}))$,  we have
$$f(\varphi_{p})\cdot y^{u}=\sum_{z\in S(\varphi_{p})}z\cdot y^{u}=\sum_{z\in S(\varphi_{p})}z^{u}\cdot y^{u}=f(\varphi_{p})\cdot y $$
because we know that $z^{u}\cdot y^{u}=z\cdot y$. Therefore we deduce that $\alpha (\varphi)$ may be written as:
$$\alpha (\varphi)=-2\sum_{1\leq j\leq  r}(f(\varphi_{p})\cdot y_{j})A_{j}(\varphi_{p})$$
with
$$A_{j}(\varphi_{p})=\sum_{1\leq i\leq q}\sum _{u\in V_{p}/V_{p,j}}
\frac{n(\varphi_{p}^{a_{i}}, y_{j}^{u})}{e_{y_{j}^{u}}}
\sigma_{a_{i}}^{-1} \ .$$
Let  $l^{m_{j}}$  denote the order of $V_{p,j}$. Using the equality
$n(\varphi_{p}^{a_{i}}, y_{j}^{u})=n(\varphi_{p}^{a_{i}u^{-1}}, y_{j})$ and the fact that $e_y=e_{y^u}$, from the above we obtain that
$$l^{m_{j}}A_{j}(\varphi_{p})=
\sum_{1\leq i\leq q}\sum _{u\in V_{p}}\frac{n(\varphi_{p}^{a_{i}}, y_{j}^{u})}{e_{y_{j}^{u}}}
\sigma_{a_{i}}^{-1} =\sum_{1\leq i\leq q}\sum _{v\in V_{p}}\frac{n(\varphi_{p}^{a_{i}v}, y_{j})}{e_{y_{j}}}
\sigma_{a_{i}}^{-1}$$
and since for $v\in Gal({\bf Q}_{p}(\zeta _{l^{N}}):{\bf Q}_{p}(\varphi_{p})))$ we have $n(\varphi^v)=n(\varphi)$ we get
$$l^{m_{j}}A_{j}(\varphi_{p})=[{\bf Q}_{p}(\zeta_{l^{N}}):{\bf Q}_{p}(\varphi_{p})](\sum_{1\leq i\leq q}\sum _{v\in H_{p}}\frac{n(\varphi_{p}^{a_{i}v}, y_{j})}{e_{y_{j}}}
\sigma_{a_{i}}^{-1}\ .$$
 Since the integer $[{\bf Q}_{p}(\zeta _{l^{N}}):{\bf Q}_{p}(\varphi_{p})]$ is a power  of $l$ we deduce that for each $j$ there exists
  $n_{j}\in {\bf Z}$ such that

$$A_{j}(\varphi_{p})=l^{n_{j}}\sum_{1\leq i\leq q}\sum _{v\in H_{p}}\frac{n(\varphi_{p}^{a_{i}v}, y_{j})}{e_{y_{j}}}
\sigma_{a_{i}}^{-1}\ .$$
  Let $l^{s}$ denote the order of $\varphi$. Since  for any $v\in H_{p}$, we know that $\mathcal{P}^{\sigma_{v}}=\mathcal{P}$,
we can  deduce from the previous equality and the definition of $\alpha(\varphi)$ that there exists a positive integer $m$ and for any $j, 1\leq j\leq q$,
an integer $a_{j}$ such that
$$c(h_{p}(\varphi))^{l^{m}}= \mathcal{P}^{l^m\alpha(\varphi)}$$
with
$${l^m\alpha(\varphi)}=-2\sum_{1\leq j\leq r}a_{j}b_{j}(\varphi)\ \quad \mbox{and}\ \quad b_{j}(\varphi)=\sum_{1\leq a \leq l^{s},
 (a,l)=1}
\frac {n(\varphi_{p}^{a}, y_{j})}{e_{y_{j}}}\sigma_{a}^{-1}\ .$$
We now wish to study the sums $b_{j}(\varphi)$. Let $y$ be one of the $\{y_{j}\}$ and let
$$b(\varphi)=\sum_{1\leq a \leq l^{s}, (a,l)=1}\frac {n(\varphi_{p}^{a}, y)}{e_{y}}\sigma_{a}^{-1}\ .$$
We write $e_{y}=l^{n}, \varphi_{p}\mid I_{y}=\psi _{y}^{ul^{m}}$ with $(u,l)=1$ and $0<ul^{m}<l^{n}$. It follows that the order of
$\varphi_{p}\mid I_{y}$ is equal to $l^{n-m}$; if we put $t=(n-m)$  and recall that we denote the order of $\varphi_{p}$  by $l^{s}$, then we have $t\leq s$.
One checks easily that
$$\frac{n(\varphi_{p}^{a}, e_{y})}{e_{y}}=\{\frac{au}{l^{t}}\}.$$ We have the equality of sets of integers:
$$\{c+kl^{t}, 1\leq c<l^{t}, (c, l)=1, 0\leq k<l^{s-t}\}=\{1\leq a< l^{s}, (a, l)=1\}\ .$$
This implies that $b(\varphi)$ may be rewritten:
$$b(\varphi)=\sum_{1\leq c<l^{t}, (c, l)=1}\{\frac {uc}{l^{t}}\}\sum_{0\leq k<l^{s-t}}\sigma_{c+kl^{t}}^{-1}\ .$$
If we let  $M_{y}$ denote the field ${\bf Q}(\zeta_{l^{t}})$, then it follows  that
\begin{eq}
b(\varphi)=\sigma_{u}\theta (M_{y})Tr _{{\bf Q}(\varphi)/M_{y}}
\end{eq}
where $\theta(M_{y})$ is the Stickelberger element of the field $M_{y}$. We then deduce from Stickelberger's Theorem that for any $j$ there exists an integer $N_{j}$
such that $l^{N_{j}}b_{j}(\varphi)$ annihilates $\mathcal{P}$ and therefore that there exists a power of $l$ which annihilates
$c(h_{p}(\varphi))$, as required for the first part of the lemma.
\bigskip

In the special case where $G$ has order $l$ we note that in the above $m_{j}=n_{j}=0$ for any $j$ and that
${\bf Q}_{p}(\zeta_{l^{N}})={\bf Q}_{p}(\varphi_{p})$. Therefore for any $j$ there exists $r_{j}$ such that
$A_{j}(\varphi_{p})=\sigma_{r_{j}}\theta$, where $\theta$ is the Stickelberger element of ${\bf Q}(\zeta_{l})/
{\bf Q}$. We therefore conclude that $c(h_{p}(\varphi))$ is principal. The lemma is now proved.
\end{Proof}
\medskip

We now can complete the proof of Theorem 1.9.  It follows from Lemmas 5.1 and 5.2  that for any $p \in \Sigma $ there exists an integer $d$
such that $t(g_{p})^{l^{d}}$ belongs to $D({\bf Z}[G])$. Moreover one can choose  $d$  equal to $0$ when $G$ is of order $l$. This
proves Theorem 1.9. $i) $ and $iii) $,  since the group $D({\bf Z}[G])$ is itself an $l$-group which is trivial when $G$ is of order $l$.
 Let $\mathcal{M}$ denote the unique maximal order of the group algebra $\bf{Q}[G]$. When $l$ is regular,  the order of the class group
$\mbox{Cl}(\mathcal{M})$ is coprime with  $l$. Therefore,  in this case,   for any $p\in \Sigma$, the class  $t(g_{p})$ is
itself an element of  $D({\bf Z}[G])$.  This proves Theorem 1.9 $(ii)$.

\bigskip
\section {Some non-trivial Euler characteristics }

      This section contains the proof of Theorem 1.10.  Our  aim here is to use the theory of modular curves to provide some detailed  examples
       of the above general results.

      \subsection {A modular example}
      We let $p$ be a prime number with $p\equiv 1$ mod $24$ and  $X_{1}=X_{1}(p)$ be the model over  $\mbox {\Spec}({\bf Z})$ of the modular curve
      associated to the congruence subgroup $\Gamma _{1}(p)$ as considered in  Section 4 of [CPT].
       The group $\Gamma /\{\pm 1\}$ acts faithfully on $X_{1}$. Let $l$ be a prime divisor of $(p-1)$ with $l>3$ and let
      $H$ be the subgroup  of $\Gamma $ of index $l$. Since $p\equiv 1$ mod $24$ we know that $6$ divides the order of $H$. We consider   the quotients $Y=X_{1}/ \Gamma $ and $X=X_{1}/H$ of $X_{1}$.
      The scheme $X$ is projective and flat over $\mbox {\Spec}({\bf Z})$ and is acted on tamely by the cyclic group $G=\Gamma /H$ of order $l$.
      The morphism $\pi: X\rightarrow Y$ is a cover which fulfils all the hypotheses required in this paper (see  Theorems 4.2  and
       4.3 in [CPT]).  It follows from the
      work of various authors that $X[1/p]\rightarrow Y[1/p]$ is a $G$-torsor which implies that  the set $\Sigma $ defined previously
       is reduced to $\{p\}$. Moreover the special fiber over $p$ is reduced with normal crossings and   has two irreducible components
       $D_{0}$ and $D_{\infty}$ both isomorphic to ${\bf P}^{1}_{{\bf F}_{p}}$. If we let $y_{0}$ (resp. $y_{\infty}$) be the generic
       point of $D_{0}$ (resp. $D_{\infty}$),  then the intersection number $y_{0}.y_{\infty}$ is equal to $(p-1)/12$. Therefore we deduce that
       $y_{0}^{2}=y_{\infty}^{2}=(1-p)/12$ (see for instance  Proposition 1.21 in Chapter 9 of [L]). Finally, it follows once again from Theorem 4.3 in [CPT]
        that $\pi$ is totally ramified over $y_{0}$ and non-ramified over $y_{\infty}$. If $x_{0}$ denotes the codimension one
       point of $X$ above $y_{0}$, then  we denote by $\psi _{0}$ the character of $G=I_{x_{0}}$ defining the action of
      $ I_{x_{0}}$  on the cotangent space at $x_{0}$.
      \medskip

      \noindent Before describing our result we need to introduce a small amount of further notation.
       We let $\mathcal{P}$ be the prime ideal of the cyclotomic field ${\bf Q}(\zeta _{l})$ defined by
      the chosen embedding $j_{p}: {\bf Q^{c}}\rightarrow {\bf Q}^{c}_{p} $. For any abelian character
      $\theta $ of $G$ we denote by $\theta _{p}$ the ${\bf Q}^{c}_{p}$-character of $G$ obtained by composing $\theta $
      with $j_{p}$. We denote by $\varphi $ the non-trivial character of $G$ such that $\varphi_{p}=\psi _{0}$. The character
      $\varphi$ induces a group homomorphism $\lambda \mapsto \varphi (\lambda )$ from  the group of ideles of ${\bf Q}[G]$
      to the group of ideles of ${\bf Q}(\zeta_{l} )$. We denote by $c(\varphi (\lambda ))$ the content of this idele.  For any number field
      $L$ and any fractional ideal $I$  of $L$ we denote by $[I]$ the class of $I$ in the class group $\mbox {Cl}(O_L)$ of $O_{L}$.
      Let $x\in \mbox{Cl}({\bf Z}[G])$ have representative $\lambda$ in the
      group of ideles of ${\bf Q}[G]$. Since $G$ is a group of  prime order $l$, the map $x\mapsto [c(\varphi (\lambda))]$,  induces a group isomorphism from $\mbox{Cl}({\bf Z}[G])$ into $\mbox{Cl}({\bf Z}[\zeta _{l}])$ that
      we also denote by $\varphi$.
      \medskip

      \noindent For $\mathcal{F}$ equal to either $\mathcal{O}_{X}$ or $\omega _{X/Y}^{1/2}$ and for any
      abelian  character $\theta $ of
      $G$ we consider the rational number $T_{p}(\mathcal{F},  \theta  )$:
      \begin{eq}
      T_{p}(\mathcal{F},  \theta  )=r_{p}(\mathcal{F},  \theta)^{2}+c_{1}(\omega _{Y_{p}})\cdot r_{p}(\mathcal{F},\theta)\ .
      \end{eq}
      We know that $l$ divides $y_0\cdot y_{\infty}$; it follows from Proposition 4.1 that $T _{p}(\mathcal{O}_{X},  \theta )-T _{p}(\omega _{X/Y}^{1/2}, \theta )$
      is an integer and from 3.13 in  [CPT] we know that the same is true for  $lT _{p}(\mathcal{O}_{X},  \theta )$. Let us first consider the class
      $$V_{X/Y}=2\chi^{P}(\omega _{X/Y}^{1/2})-2\chi^{P}(\mathcal{O}_{X})\ .$$
      By Theorem 1.7  we conclude that
      $$ \varphi(V_{X/Y})=[\mathcal{P}]^{\sum_{1\leq a<l}
      ({T_{p}(\mathcal{O}_{X},  \psi _{0}^{a} )-T _{p}(\omega _{X/Y}^{1/2},  \psi_{0}^{a} )})
      \sigma_{a}^{-1}} \ \ \ \ .$$
      The description of  $\varphi (l\chi(\mathcal{O}_{X})$ follows from Theorem 3.13 in [CPT].
       Since $lT _{p}(\mathcal{O}_{X},   \theta )$ is an integer, for any abelian character  $\theta$ we obtain that

       $$\varphi(2l\chi (\mathcal{O}_{X}))=
       [\mathcal{P}]^{-l\sum_{1\leq a<l}T _{p}(\mathcal{O}_{X}, \psi _{0}^{a} )\sigma_{a}^{-1}}  . $$
       We deduce from  the previous equalities  that
      $$ \varphi(2l\chi (\omega_{X/Y})^{1/2})=
       [\mathcal{P}]^{-l\sum_{1\leq a<l}T _{p}(\omega _{X/Y}^{1/2}, \psi _{0}^{a} )\sigma_{a}^{-1}}\ .$$

       \subsection {Proof of Theorem 1.10}
        From 3.5 we know that
        $r_{p}(\mathcal{O}_{X},  \psi _{0}^{a})=\frac {a}{l}y_{0}$.
        On the other hand  from (3.15) in Section 3  in [CPT] we know that
        \begin{eq}
        c_{1}(\omega_{Y_{p}})\cdot y_{0}=-y_{0}^{2}-2\chi(y_{0}, \mathcal{O}_{y_{0}})=-y_{0}^{2}-2.
        \end{eq}
        Piecing this together with 6.1 we see that we have shown that
         $$-lT_{p}(\mathcal{O}_{X}\mid \psi _{0}^{a})=\frac {p-1}{12l}a^{2}+(2- \frac {p-1}{12})a$$
         and so
         $$\varphi (2l\chi(\mathcal{O}_{X}))=
         [\mathcal{P}]^{\frac{p-1}{12l}\sum_{1\leq a<l}a^{2}\sigma_{a}^{-1}+(2-\frac {p-1}{12})l\theta} $$
         where $\theta$ is the Stickelberger element of ${\bf Q}(\zeta_{l})$.
         Since we know from Stickelberger's theorem that $l\theta$ annihilates the class group of ${\bf Q}(\zeta _{l})$, we deduce that
         $$\varphi (2l\chi(\mathcal{O}_{X}))=
         [\mathcal{P}]^{\frac{p-1}{12l}\sum_{1\leq a<l}a^{2}\sigma_{a}^{-1}}. $$
         We apply  the results of Examples in Section 3 and consider separately the cases when $a<l/2$ and when $a>l/2$.

         If $a<l/2$, then we obtain that
         $$r_{p}(\omega _{X/Y}^{1/2},  \psi _{0}^{a})=r_{p}(\mathcal{O}_{X}, \psi _{0}^{a})$$
         and so   $$T _{p}(\mathcal{O}_{X}, \psi _{0}^{a} )-T _{p}(\omega _{X/Y}^{1/2},  \psi_{0}^{a} )=0.$$
         If now  $a>l/2$, then
         $$r_{p}(\omega _{X/Y}^{1/2},  \psi _{0}^{a})=r_{p}(\mathcal{O}_{X},  \psi _{0}^{a})-y_{0}$$
         and so
         $$T _{p}(\mathcal{O}_{X}, \psi _{0}^{a} )-T _{p}(\omega _{X/Y}^{1/2},  \psi_{0}^{a} )=\frac {p-1}{6}(1-\frac {a}{l})-2. \ $$

          Therefore we have proved that
          $$\varphi(V_{X/Y})=\overline \mathcal{P}^{\frac {p-1}{6l}s_{1}-2s_{0}} \ $$
           where for any integer $i\in \{0, 1, 2\}$ we have defined   $s_{i}$ in  ${\bf Z}[G]$ by
       $$s_{i}=\sum_{1\leq a<l/2}a^{i}\sigma_{a}^{-1}\ .$$
      In order to conclude we need some elementary relationships in ${\bf Z}[G]$ that we prove in the next lemma.
         \begin {lemma} One has the following equalities:

         \noindent i) $s_{0}=(2\sigma_{-1}-\sigma_{(l-1)/2}^{-1})\theta $.

        \noindent  ii) $(1-\sigma_{-1})s_{1}=(\sigma_{(l+1)/2}^{-1}-1)(l\theta)$.

        \noindent iii) $\sum_{1\leq a<l}a^{2}\sigma_{a}^{-1}=(1+\sigma_{-1})s_{2}+l\sigma_{-1}
        (ls_{0}-2s_{1})$.
        \end {lemma}
        \begin{Proof} The proof of $iii)$ is immediate. We deduce $ii)$ from $ i)$ via the  relationship:
       $$l\theta =(1-\sigma_{-1})s_{1}+l\sigma_{-1}s_{0}\ .$$
       Therefore it suffices to prove $i)$. We first observe that
      $$2\sigma_{2}^{-1}(l\theta)=\sum_{1\leq a<l}(2a\sigma_{2a}^{-1})$$
      can be decomposed as the sum
      $$2\sigma_{2}^{-1}(l\theta)=\sum_{1\leq a<l/2}2a\sigma_{2a}^{-1}+ \sum_{l/2<a<l}(2a-l)\sigma_{2a-l}^{-1}
      +l\sum_{l/2<a<l}\sigma_{2a-l}^{-1}\ .$$
      We now observe that the sum of the first two terms is precisely $l\theta $. Moreover, using the equality
      $\sigma^{-1} _{2a-l}=\sigma_{\frac {l-1}{2}}\sigma_{l-a}^{-1}$, we express the last term of the sum as
      $\sigma_{\frac {l-1}{2}}s_{0}$. Therefore we have proved that
      $$2\sigma_{2}^{-1}\theta =\theta+\sigma_{\frac {l-1}{2}}s_{0}\ $$
      and  $i)$ follows.
      \end{Proof}
      \medskip

      From the lemma and Stickelberger's theorem we deduce  that $s_{0}$ annihilates the class group of
      ${\bf Z}[\zeta _{l}]$. Therefore we obtain that
      $$\varphi (V_{X/Y})=[\mathcal{P}]^{\frac {p-1}{6l}s_{1}}\ \mbox {and}
      \ \ \varphi (2l(\chi(\mathcal{O}_{X})))=
      [\mathcal{P}\overline \mathcal{P}]^{\frac {p-1}{12l}s_{2}}[\overline \mathcal{P}]^{-\frac {p-1}{6}s_{1}}\ .$$
      It follows also from part $(ii)$ of the lemma  that $[\mathcal{P}]^{s_{1}}=[\overline \mathcal{P}]^{s_{1}}$ and thus  that
      $[\mathcal{P}]^{2s_{1}}=[\mathcal{P}\overline \mathcal{P}]^{s_{1}}$. The formulas of Theorem 1.10  follow immediately.

      \subsection {Non existence of NIB}

      We conclude this section  by   giving examples where our invariants are non-trivial. Our strategy is to evaluate the norm of these  classes in
      the class group of the quadratic subfield $k$ of ${\bf Q}(\zeta _{l})$.  We will denote  the norm from ${\bf Q}(\zeta _{l})$ to $k$ by $N$.
       Since ${\bf Q}(\zeta _{l})/k$ contains no unramified subextension $F/k$ with $F\not=k$,  $N$ induces a surjective
       group homomorphism from $\mbox{Cl}({\bf Z}[\zeta _{l}])$ onto $\mbox {Cl}(\mathcal O_k)$ (see for instance Theorem 10.1 of [W]).  When $l\equiv 3$ mod $4$ the field $k$ is quadratic imaginary.
       It therefore follows from Theorem 1.10.i  that all the  classes belong to the kernel of $N$. This implies a certain restriction on
       their orders.
       \medskip

       We next consider the case  where $l\equiv  1$ mod $4$, so that  $k={\bf Q}(\sqrt l)$. The Galois group of
       ${\bf Q}(\zeta _{l})$ over $k$ consists of the set $\{\sigma _{a}, 1\leq a<l\}$ such that $a$ is a square mod $l$. We let
       $A$ (resp. $B$) denote the set of $a, 1\leq a<l/2$,  such that $a $ is (resp. is not ) a square mod $l$. Since $-1$ is
       a square, we see immediately
      that $\rm{Card} {(A)}=\rm{Card} {(B)}$.
      For $i\in \{1, 2\}$ we set $$t_{i}=\sum_{a\in A}a^{i}-\sum_{b\in B}b^{i}.$$ By
      taking the norm of both sides of the equalities in Theorem 1.10, and writing  $\beta =N(\mathcal{P})$, we obtain that
      $$N(\varphi(V_{X/Y}))=[\beta]^{\frac{p-1}{6l}t_{1}}\ ,\quad
      N(\varphi( 2l\chi (\omega _{X/Y}^{1/2})))=[\beta]^{\frac{p-1}{6l}t_{2}}$$
      {and}
    $$N(\varphi( 2l\chi (\mathcal{O}_{X})))=[\beta]^{\frac{p-1}{6l}(t_{2}-lt_{1})}.$$
    \bigskip

    \noindent To conclude we consider the case  $l=401$.  The class number of ${\bf Q}(\sqrt l)$ is $5$. The integers $t_{1}$ and $t_{2}$
    are independent of $p$. By computation we obtain that $t_{1}\equiv 4$ mod $5$, $t_{2}\equiv 3$ mod $5$ and
    $t_{2}-lt_{1}\equiv  4$ mod $5$. Therefore for any $p\equiv 1$ mod $24l$ with $p\not\equiv 1$ mod $5$ with the property that
    $\beta $ is not principal in $k$, we obtain three non trivial classes. The smallest prime satisfying these properties is
    $p=182857$. This example therefore  provides us with a tame cover of surfaces $\pi: X\rightarrow Y$ where
    $2\chi (\mathcal{O}_{X})$, $2\chi(\omega _{X/Y})$ and $2\chi(\omega _{X/Y}^{1/2})$ are all non trivial,
    where $2\chi (\mathcal{O}_{X})=2\chi(\omega _{X/Y})$ but where
     $2\chi(\mathcal{O}_{X/Y})\not\neq 2\chi(\omega _{X/Y}^{1/2})$.
   
philippe.cassou-nogues@math.u-bordeaux1.fr
\\
\bigskip
martin.taylor@manchester.ac.uk

\bigskip

     \begin{thebibliography}{C-P-S 123} \addcontentsline{toc}{section}{Bibliography}

\bibitem [C1] {C1} T. Chinburg,
{\it Galois structure of the de Rham cohomology of tame covers of schemes},  Ann. of Math. 139 (1994), 443-490.

\bibitem [C2] {C2} T. Chinburg,
{\it Galois module structure of de Rham cohomology}, J. de Th\'eorie des Nombres de Bordeaux 4 (1991), 1-18.

\bibitem [CE] {CE} T. Chinburg, B. Erez,
{\it Equivariant Euler-Poincar\'e characteristics and tameness},  Ast\'erisque 209 (1992) 179-194.

\bibitem [CEPT1] {CEPT1} T. Chinburg, B. Erez, G. Pappas, M.J. Taylor,
{\it $\varepsilon $-constants and the Galois structure of de Rham cohomology}, Ann. of Math. 146 (1997), 411-473.

\bibitem [CEPT2] {CEPT2} T. Chinburg, B. Erez, G. Pappas, M.J.  Taylor,
{\it Riemann-Roch type theorems for arithmetic schemes with a finite group action}, J. Reine Angew. Math. 489 (1997), 151-187.

\bibitem [CEPT3] {CEPT3} T. Chinburg, B. Erez, G. Pappas, M.J.  Taylor,
{\it Tame actions of group schemes: integrals amd slices}, Duke Math. J. 82 (1996), 269-308.

\bibitem [CPT] {CPT} T. Chinburg, G. Pappas, M.J.  Taylor,
{\it Cubic structures, equivariant Euler characteristics and lattices of modular forms}, Ann. of Math. 170 (2009), 561-608.

\bibitem [ET] {ET} B. Erez, M.J. Taylor,
{\it Hermitian modules in Galois extensions of number fields and Adams operations}, Ann. of Math. 135 (1992), no. 2, 271-296.

\bibitem [F1]  {F}  A. Fr\"ohlich,
{\it Galois module structure of algebraic integers},  Ergebnisse der Mathematik und ihrer Grenzgebiete (3) 1,
Springer-Verlag  (1983).

\bibitem [F2]  {F1} A. Fr\"ohlich,
{\it Arithmetic and Galois module structure for tame extensions}, J. Crelle 286/287 (1976), 380-440.

\bibitem [Fu] {Fu} W. Fulton,
{\it Intersection theory}, 2nd edition, Ergebnisse der Mathematik und ihrer Grenzgebiete  3. Folge, Springer-Verlag (1998).
\bibitem [H] {H} R. Hartshorne,
{\it Algebraic Geometry}, Grad. Texts  Math. 52, Springer-Verlag (1977).

\bibitem [RD]  {RD} R. Hartshorne,
{\it Residues and Duality },  Lect.  Notes in Math. 20,  Springer-Verlag (1966).



\bibitem [L] {L} S. Lang,
{\it Introduction to Arakelov theory}, Springer-Verlag (1988).


\bibitem [P]  {P} G. Pappas,
{\it Galois modules and the Theorem of the Cube},  Invent. Math. 133 (1998), 193-225.

\bibitem [T1] {T1} M.J. Taylor,
{\it On the self duality of a ring of integers as a Galois module}, Invent. Math. 46 (1978), 173-177.

\bibitem [T2] {T2} M.J. Taylor,
{\it Class groups of group rings},  London Math. Soc. Lect.  Notes 91, C. U. P. (1983).

\bibitem [T3] {T1} M.J. Taylor,
{\it On Fr\"ohlich's conjecture for rings of integers of tame extensions}, Invent. Math. 63 (1981), 41-79.


\bibitem [W] {W} L. Washington,
{\it Introduction to cyclotomic fields}, Graduate Texts in Mathematics, Springer Verlag, 1980.







\end{thebibliography}
\end{document}